\documentclass{conm-p-l}
\usepackage{amssymb,latexsym,amscd}
\newcommand{\Mdef}[2]{\newcommand{#1}{\relax \ifmmode #2 \else $#2$\fi}}

\newcommand{\cok}{\mathrm{cok}}

\newcommand{\sm }{\wedge}

\newcommand{\tensor}{\otimes}

\newcommand{\Hom}{\mathrm{Hom}}

\newcommand{\Ext}{\mathrm{Ext}}
\Mdef{\bhom}{\mathbf{\hat{H}om}}

\Mdef{\Mod}{\mathrm{mod}}

\newcommand{\st}{\; | \;}

\newtheorem{thm}{Theorem}[section]
\newtheorem{lemma}[thm]{Lemma}
\newtheorem{prop}[thm]{Proposition}
\newtheorem{cor}[thm]{Corollary}

\theoremstyle{definition}

\newtheorem{defn}[thm]{Definition}

\newtheorem{context}[thm]{Context}

\newtheorem{example}[thm]{Example}
\newtheorem{examples}[thm]{Examples}
\newtheorem{remark}[thm]{Remark}

\newcommand{\qqed}{\qed \\[1ex]}
\renewenvironment{proof}[1][\hspace*{-.8ex}]{\noindent {\bf Proof #1:\;}}{\qqed}

\Mdef{\PH} {\Phi^H}
\Mdef{\PK} {\Phi^K}
\Mdef{\PL} {\Phi^L}
\Mdef{\PT} {\Phi^{\T}}

\Mdef{\ef}{E{\cF}_+}
\Mdef{\etf}{\tilde{E}{\cF}}
\Mdef{\eg}{E{G}_+}
\Mdef{\etg}{\tilde{E}{G}}
\Mdef{\tf}{\T / \! \!  / {\cF}_+}

\Mdef{\infl}{\mathrm{inf}}
\Mdef{\defl}{\mathrm{def}}
\Mdef{\res}{\mathrm{res}}
\Mdef{\ind}{\mathrm{inf}}

\Mdef{\univ}{\mathcal{U}}

\Mdef{\Fp}{\mathbb{F}_p}
\Mdef{\Zpinfty}{\Z /p^{\infty}}
\Mdef{\Zpadic}{\Z_p^{\wedge}}

\newcommand{\bi}{\begin{itemize}}
\newcommand{\be}{\begin{enumerate}}
\newcommand{\bc}{\begin{center}}
\newcommand{\bd}{\begin{description}}
\newcommand{\ei}{\end{itemize}}
\newcommand{\ee}{\end{enumerate}}
\newcommand{\ec}{\end{center}}
\newcommand{\ed}{\end{description}}

\newcommand{\adjunction}[4]{
\diagram
#1:#2 \rrto<0.7ex> &&
#3  \llto<0.7ex> :#4 
\enddiagram}
\newcommand{\lra}{\longrightarrow}
\newcommand{\lla}{\longleftarrow}

\newcommand{\iso}{\cong}

\Mdef{\we}{\mathbf{we}}
\Mdef{\fib}{\mathbf{fib}}
\Mdef{\cof}{\mathbf{cof}}
\Mdef{\BI}{\mathcal{BI}}

\newcommand{\ilim}{\mathop{ \mathop{\mathrm{lim}} \limits_\leftarrow} \nolimits}
\newcommand{\colim}{\mathop{  \mathop{\mathrm {lim}} \limits_\rightarrow} \nolimits}
\newcommand{\holim}{\mathop{ \mathop{\mathrm {holim}} \limits_\leftarrow} \nolimits}
\newcommand{\hocolim}{\mathop{  \mathop{\mathrm {holim}}\limits_\rightarrow} \nolimits}

\Mdef{\A}{\mathbb{A}}
\Mdef{\B}{\mathbb{B}}
\Mdef{\C}{\mathbb{C}}
\Mdef{\D}{\mathbb{D}}
\Mdef{\E}{\mathbb{E}}
\Mdef{\T}{\mathbb{T}}
\Mdef{\F}{\mathbb{F}}
\Mdef{\G}{\mathbb{G}}
\Mdef{\I}{\mathbb{I}}
\Mdef{\N}{\mathbb{N}}
\Mdef{\Q}{\mathbb{Q}}
\Mdef{\R}{\mathbb{R}}
\Mdef{\bbS}{\mathbb{S}}
\Mdef{\Z}{\mathbb{Z}}

\Mdef{\bA}{\mathbb{A}}
\Mdef{\bB}{\mathbb{B}}
\Mdef{\bC}{\mathbb{C}}
\Mdef{\bD}{\mathbb{D}}
\Mdef{\bE}{\mathbb{E}}
\Mdef{\bF}{\mathbb{F}}
\Mdef{\bG}{\mathbb{G}}
\Mdef{\bH}{\mathbb{H}}
\Mdef{\bI}{\mathbb{I}}
\Mdef{\bJ}{\mathbb{J}}
\Mdef{\bK}{\mathbb{K}}
\Mdef{\bL}{\mathbb{L}}
\Mdef{\bM}{\mathbb{M}}
\Mdef{\bN}{\mathbb{N}}
\Mdef{\bO}{\mathbb{O}}
\Mdef{\bP}{\mathbb{P}}
\Mdef{\bQ}{\mathbb{Q}}
\Mdef{\bR}{\mathbb{R}}
\Mdef{\bS}{\mathbb{S}}
\Mdef{\bT}{\mathbb{T}}
\Mdef{\bU}{\mathbb{U}}
\Mdef{\bV}{\mathbb{V}}
\Mdef{\bW}{\mathbb{W}}
\Mdef{\bX}{\mathbb{X}}
\Mdef{\bY}{\mathbb{Y}}
\Mdef{\bZ}{\mathbb{Z}}

\Mdef{\cA}{\mathcal{A}}
\Mdef{\cB}{\mathcal{B}}
\Mdef{\cC}{\mathcal{C}}
\Mdef{\mcD}{\mathcal{D}} % Something funny about \cD.
\Mdef{\cE}{\mathcal{E}}
\Mdef{\cF}{\mathcal{F}}
\Mdef{\cG}{\mathcal{G}}
\Mdef{\mcH}{\mathcal{H}} % There's something funny about \cH: it 
\Mdef{\cI}{\mathcal{I}}
\Mdef{\cJ}{\mathcal{J}}
\Mdef{\cK}{\mathcal{K}}
\Mdef{\mcL}{\mathcal{L}}% There's something funny about \cL: it 
\Mdef{\cM}{\mathcal{M}}
\Mdef{\cN}{\mathcal{N}}
\Mdef{\cO}{\mathcal{O}}
\Mdef{\cP}{\mathcal{P}}
\Mdef{\cQ}{\mathcal{Q}}
\Mdef{\mcR}{\mathcal{R}}% There's something funny about \cR: it 
\Mdef{\cS}{\mathcal{S}}
\Mdef{\cT}{\mathcal{T}}
\Mdef{\cU}{\mathcal{U}}
\Mdef{\cV}{\mathcal{V}}
\Mdef{\cW}{\mathcal{W}}
\Mdef{\cX}{\mathcal{X}}
\Mdef{\cY}{\mathcal{Y}}
\Mdef{\cZ}{\mathcal{Z}}

\Mdef{\tA}{\tilde{A}}
\Mdef{\tB}{\tilde{B}}
\Mdef{\tC}{\tilde{C}}
\Mdef{\tE}{\tilde{E}}
\Mdef{\tH}{\tilde{H}}
\Mdef{\tK}{\tilde{K}}
\Mdef{\tL}{\tilde{L}}
\Mdef{\tM}{\tilde{M}}
\Mdef{\tN}{\tilde{N}}
\Mdef{\tP}{\tilde{P}}

\Mdef{\Ab}{\overline{A}}
\Mdef{\Bb}{\overline{B}}
\Mdef{\Cb}{\overline{C}}
\Mdef{\Db}{\overline{D}}
\Mdef{\Hb}{\overline{H}}
\Mdef{\Ib}{\overline{I}}
\Mdef{\Kb}{\overline{K}}
\Mdef{\Lb}{\overline{L}}
\Mdef{\Mb}{\overline{M}}
\Mdef{\Nb}{\overline{N}}
\Mdef{\Qb}{\overline{Q}}
\Mdef{\Tb}{\overline{T}}

\Mdef{\qb}{\overline{q}}
\Mdef{\rb}{\overline{r}}
\Mdef{\tb}{\overline{t}}
\Mdef{\ub}{\overline{u}}
\Mdef{\vb}{\overline{v}}
\Mdef{\hc}{\hat{c}}
\Mdef{\he}{\hat{e}}
\Mdef{\hf}{\hat{f}}
\Mdef{\hA}{\hat{A}}
\Mdef{\hH}{\hat{H}}
\Mdef{\hJ}{\hat{J}}
\Mdef{\hM}{\hat{M}}
\Mdef{\hP}{\hat{P}}
\Mdef{\hQ}{\hat{Q}}

\Mdef{\bolda}{\mathbf{a}}
\Mdef{\boldb}{\mathbf{b}}
\Mdef{\boldD}{\mathbf{D}}

\Mdef{\fm}{\frak{m}}

\Mdef{\eps}{\epsilon}

\input{xypic}

\newcommand{\id}{\mathrm{id}}
\newcommand{\tand}{\mbox{ and }}
\newcommand{\sE}{\cE}
\newcommand{\epz}{\varepsilon}
\newcommand{\sset}[1]{#1}
\newcommand{\SI}{\Sigma}
\newcommand{\htp}{\simeq}

\newcommand{\darrow}{\longrightarrow}

\newcommand{\AL}{\mbox{\boldmath $\alpha$}} \newcommand{\FK}[1]{K^{\bullet}( #1)} \newcommand{\PPK}[1]{PK^{\bullet}( #1)} \newcommand{\FC}[1]{\check{C}^{\bullet}( #1)} \newcommand{\PC}[1]{P\check{C}^{\bullet}( #1)}

\newcommand{\CJI}{\check{C}(I)}
\newcommand{\CIM}{M[I^{-1}]}

\newcommand{\hR}{\boldR}
\newcommand{\boldR}{R}%{\mathbf{R}}
\newcommand{\boldk}{k}%{\mathbf{k}}
\newcommand{\HombR}{\Hom_{\boldR}}
\newcommand{\HomR}{\Hom_{R}}

\newcommand{\kvee}{k^{\vee}}

\newcommand{\tensorcE}{\tensor_\cE}
\newcommand{\tensorR}{\tensor_R}
\newcommand{\modcE}{\mathrm{mod}\!-\!\cE}
\newcommand{\bRmod}{\boldR\!-\!\mathrm{mod}}
\gdef\overto#1{{\buildrel{#1}\over\longrightarrow}}

\newcommand{\cell}{\mathrm{Cell}}

\newcommand{\Gk}{\Gamma_k}
\newcommand{\Lk}{\Lambda^k}
\newcommand{\Rmod}{\mbox{$R$-mod}}

\newcommand{\End}{\mathrm{End}} 
\newcommand{\ExtR}{\mathrm{Ext}_R}

\newcommand{\HomE}{\mathrm{Hom}_{\cE}}

\newcommand{\Homk}{\mathrm{Hom}_k}

\newcommand{\tensorE}{\otimes_{\cE}}

\newcommand{\Ep}{E'}
\newcommand{\ksharp}{k^{\#}}
\newcommand{\depth}{\mathrm{depth}}
\newcommand{\AR}{R}
\newcommand{\kbar}{\overline{k}}
\newcommand{\tR}{\tilde{R}}
\newcommand{\tX}{\tilde{X}}
\newcommand{\HomtR}{\mathrm{Hom}_{\tR}}

\title{First steps in brave new commutative algebra.}
\author{J.P.C. Greenlees}
\address{The University of Sheffield, Sheffield, S3 7RH\, UK} 
\email{j.greenlees@sheffield.ac.uk}
\begin{document}

\baselineskip=15pt

%    Blank box placeholder for figures (to avoid requiring any
%    particular graphics capabilities for printing this document).
\newcommand{\blankbox}[2]{%
  \parbox{\columnwidth}{\centering
%    Set fboxsep to 0 so that the actual size of the box will match the
%    given measurements more closely.
    \setlength{\fboxsep}{0pt}%
    \fbox{\raisebox{0pt}[#2]{\hspace{#1}}}%
  }%
}

%\thanks{}
\date{}

\maketitle

\tableofcontents

\section{Introduction}
To begin with, it is helpful to  explain the title. The phrase `brave new
rings' was coined by F.Waldhausen, presumably to capture
both an optimism about the possibilities of generalizing rings to ring
spectra, and a proper awareness of the risk that the new step in 
abstraction would take the subject dangerously far from its justification in
examples. This article is about doing commutative algebra with ring spectra, 
and is extensively illustrated by examples.
The `First steps' of the title suggests its design is controlled by
expository imperatives, but in fact the first steps also follow 
one particular line of historical development. From either point of view, 
this gives a coherent story.

From the expository point of view,  we could explain the choice as follows:
 completion and
localization are very basic constructions in commutative algebra, and
their formal nature makes them natural things to investigate in derived
categories. Part 1 describes a naive approach through elements, and 
Part 2 moves on to a more conceptual approach, both illustrated by examples.

From the historical point of view, one route begins with the Atiyah-Segal
completion theorem about the  $K$-theory  of  classifying spaces
 \cite{AtiyahSegal}. For a finite group $G$ this is the statement that 
$$K^0(BG)=K^0_G(EG)=K^0_G(pt)_J^{\wedge}=R(G)_J^{\wedge}$$
where $R(G)$ is the complex representation ring and $J$ is its augmentation 
ideal. It is natural to try to lift
the algebraic statements to structural statements at the level of derived
categories of spectra. In this context, the philosophy is that the theorem 
states that an algebraic and a geometric completion agree. The algebra this 
uncovers then suggests a new approach to the Atiyah-Segal completion theorem
\cite{KEG}, by starting with a version for homology. It turns out that this 
phenomenon occurs in many contexts and in some of the more algebraic ones
the role of  Morita theory begins to emerge. This leads on to a more 
conceptual approach, which, like Gorenstein rings \cite{Bass}, 
turns out to occur ubiquitously.

This article is based on the lectures I gave at the meeting in
Chicago in Summer 2004. The points of view described have developed in 
joint work with Benson, Dwyer, Iyengar  and May \cite{DGI,kappaI,GM,GML};
I am grateful to them all for numerous conversations and insights. 

\part{Localization and completion for ideals.}

Localization and completion are basic constructions in commutative algebra, 
and it is useful to be able to use them in stable homotopy theory. However, 
they can be viewed as just one special case of categorically similar 
constructions. Several important theorems can then be viewed as change of 
base results in this broadened context.

We begin with constructions in commutative algebra, and work towards
formulating them in a way which generalizes to other derived categories.
Consider an ideal $I$ in a commutative 
ring $\AR $, and then attempt to approximate modules by $I$-torsion 
modules. If we use {\em inverse} limits, we reach the notion of
the completion 
$$M_I^{\wedge}=\ilim_k M/I^kM$$
 of an $R$-module $M$.  
The algebraic fact that completion is not exact in the non-Noetherian case 
means that in treating homological invariants,
we are  forced  to work with the derived functors of completion. In fact 
topology suggests a homotopy invariant {\em construction} which 
 leads in turn to a means of {\em calculation} of the left derived functors
of completion (local homology); this then  models the 
the completion of spectra by construction. 

We began by considering {\em inverse} limits of torsion modules 
since completion is exact on Noetherian modules, so in that case the zeroth
 derived functor agrees with the original and the higher derived functors
are zero. On the other hand, if we attempt to 
approximate modules by {\em direct} limits of torsion 
modules we obtain
$$\Gamma_IM =\colim_k \HomR (A/I^k,M)=\{x \st I^kx=0 \mbox{ for } k >>0\}.$$
This time  the original functor is often zero (in the torsion free case for 
example), and it is the higher derived functors that are better behaved. These 
are calculated by a dual construction, which leads to Grothendieck's 
local cohomology modules.

We will begin in Section \ref{secloccoh} 
with the familiar commutative algebra, and then adapt it to 
ring spectra and study it in Sections \ref{sectoploccoh} to \ref{sec:locaway}. 
Finally we consider two specializations: that to 
$MU$-module spectra in Section \ref{sec:bordism}, when we obtain well known chromatic constructions, and
the motivating examples in equivariant topology in Section \ref{sec:compthms1}.
Related surveys are given in  \cite{Handbook1, Handbook2}.

\section{Algebraic definitions: Local and \v{C}ech cohomology and homology} \label{secloccoh}

The material in this section is based on \cite{G, GML,Tateca}. Background in commutative
algebra can be found in \cite{Mat,BrunsHerzog}.
\subsection{The functors}
Suppose to begin with that $\AR $ is a commutative Noetherian ring and that 
$I=(\alpha_1, \dots , \alpha_n)$ is an ideal in $\AR $. 
We shall be concerned especially with two naturally occurring functors on 
$\AR $-modules: the $I$-power torsion functor and the $I$-adic completion functor.

The $I$-power torsion functor $\Gamma_I$ is defined by 
$$M \longmapsto \Gamma_I (M) = \{ x \in M \mid I^kx=0
\mbox{ for  } k >> 0  \}.$$
We say that $M$ is an $I$-power torsion module if $M=\Gamma_I M$. 
It is easy to check that the functor $\Gamma_I$ is left exact.

Recall that the support of $M$ is the set of prime ideals $\wp$ of $\AR $ such that the localization $M_{\wp}$ is non-zero. 
We say that $M$ is supported over $I$ if every prime in the support of $M$ contains $I$. This is 
equivalent to the condition that $M[1/{\alpha}]=0$ for each $\alpha \in I$. It follows that $M$ is an 
$I$-power torsion module if and only if the support of $M$ lies over $I$.

The $I$-adic completion functor is defined by 
$$M \longmapsto M_I^{\sm } = \ilim _k M/I^kM.$$
The Artin-Rees lemma implies that $I$-adic completion is exact on finitely generated modules, but it is 
neither right nor left exact in general. 

Since the functors that arise in topology are exact functors on triangulated categories, 
we need to understand the algebraic functors at the level of the derived category, 
which is to say that we must understand their derived functors. The connection with topology 
comes through one particular way of calculating the derived functors $R^*\Gamma_I$ of 
$\Gamma_I$  and $L_*^I$ of $I$-adic completion. It also provides a 
connection between the two sets of derived functors and makes useful techniques available.

\subsection{The stable Koszul complex.}
We begin with a sequence $\sset{\alpha_1,\ldots,\alpha_n}$ of elements of $\AR $
and define various chain complexes. In Subsection \ref{subsec:invariance}
we explain why the chain complexes only depend on the radical of the ideal
$I=(\alpha_1, \ldots, \alpha_n)$ generated by the sequence, in Subsection 
\ref{subsec:lochomcoh} we define associated homology groups,  
and in Subsection \ref{subsec:derived} we give conceptual 
interpretations of this homology under Noetherian hypotheses. 

We begin with a single element $\alpha \in \AR $, and an integer $s \geq 0$, and
define the $s$th unstable Koszul complex by 
$$K_s^{\bullet}(\alpha )=(\alpha^s :\AR  \darrow  \AR  )$$
where the non-zero modules are in cohomological degrees 0 and 1. 
These complexes form a direct system as $s$ varies,
$$\begin{array}{rcccc}
K_1^{\bullet}(\alpha)&=&
\left( \right. \;\; \AR & \stackrel{\alpha}\lra &\AR \left. \; \; \right)\\
\downarrow \hspace*{4ex}&&=\downarrow && \downarrow \alpha\\
K_2^{\bullet}(\alpha)&=&
\left( \right.\;\;\AR & \stackrel{\alpha^2}\lra &\AR \left. \; \; \right)\\
\downarrow \hspace*{4ex}&&=\downarrow && \downarrow \alpha\\
K_3^{\bullet}(\alpha)&=&
\left( \right. \; \; \AR & \stackrel{\alpha^3}\lra &\AR \left. \; \;\right)\\
\downarrow \hspace*{4ex} &&=\downarrow && \downarrow \alpha
\end{array}$$
 and the direct limit is  the {\em flat stable} Koszul complex 
$$\FK{\alpha} = \left( \AR  \darrow  \AR [1/\alpha ] \right). $$ 

When defining local cohomology, it is usual to use the complex $\FK{\alpha}$ 
of flat modules. However, we shall need a complex of projective $\AR$-modules 
to define the dual local homology modules. Accordingly, we take a 
particularly convenient projective approximation $\PPK{\alpha}$ to $\FK{\alpha}$. 
Instead of taking the direct limit of the $K^{\bullet}_s(\alpha )$, we take their 
homotopy direct limit. This makes the translation to the  topological context
straightforward. More concretely, our model for $\PPK{\alpha}$ is
displayed as the upper row in the homology isomorphism 
$$\diagram PK^{\bullet}(\alpha)  \dto&=& \left( \right. 
\AR  \oplus \AR [x] \dto_{\langle 1,0 \rangle} 
\rrto^<(0.35){\langle 1,\alpha x-1\rangle} & & 
\AR [x] \dto^{g} \left.\right)\\ 
K^{\bullet}(\alpha) &=& \left( \right. \hspace*{3ex}\AR \hspace*{3ex} \rrto & &\AR [1/\alpha] \left. \right), \\ 
\enddiagram$$ 
where $g(x^i)=1/\alpha^i$. Like $\FK{\alpha}$, this choice of $\PPK{\alpha}$ is 
non-zero only in cohomological degrees $0$ and $1$.

The stable Koszul cochain complex for a sequence $\AL = (\alpha_1, \dots , \alpha_n ) $ 
is obtained by tensoring together the complexes for the  elements, so that 
$$ \FK{\AL} =  \FK{\alpha_1} \otimes_R \dots \otimes_R \FK{\alpha_n},$$ 
and similarly for the projective complex $\PPK{\AL}$.

\subsection{Invariance statements.}
\label{subsec:invariance}

We prove some basic properties of the stable Koszul complex. 

\begin{lemma}\label{Ksupp}
If $\beta$ is in the ideal $I=(\alpha_1 , \alpha_2, \ldots , \alpha_n )$, then 
$\FK{\AL }[1/\beta ]$ is exact.
\end{lemma}
\begin{proof} 
Since homology commutes with direct limits, it suffices to show that some power of $\beta$ 
acts as zero on the homology of 
$K^{\bullet}_s(\AL)=K^{\bullet}_s(\alpha_1)\otimes\cdots\otimes K^{\bullet}_s(\alpha_n)$. However, 
$(\alpha_i)^{s}$ annihilates $H^*(K^{\bullet}_s(\alpha_i))$, and it follows 
from the long exact sequence in homology that $(\alpha_i)^{2s}$ 
annihilates $H^*(K^{\bullet}_s(\AL))$. Writing $\beta$ as a linear combination of the $n$ elements 
$\alpha_i$, we see that $\beta^{2sn}$ is a linear combination of elements each of which is divisible 
by some $(\alpha_i)^{2s}$, and the conclusion follows. \end{proof}

Note that, by construction, we have an augmentation map
$$\epz : \FK{\AL } \darrow  \AR .$$

\begin{cor} \label{indepgens}
Up to quasi-isomorphism, the complex $\FK{\AL}$ depends only on the radical of the ideal $I$. 
\end{cor} 
\begin{proof} Since $\FK{\AL}$ is unchanged if we replace the generators by powers, it
 suffices to show that $\FK{\AL}$ only depends on the ideal generated by $\AL$.

The augmentation gives a map 
$$\FK{\AL , \beta } =\FK {\AL} \tensor_R \FK {\beta} \darrow  \FK{\AL}\tensor_R R=\FK{\AL}$$ 
which is a quasi-isomorphism if $\beta\in I$ since the lemma shows its cofibre
$\FK{\AL}[1/\beta]$ is exact. It follows that we have homology isomorphisms 
$$\FK{\AL} \lla  \FK{\AL} \otimes \FK{\AL'} \darrow  \FK{\AL'}$$ if $\AL'$ is a 
second set of generators for $I$. 

\end{proof}

In view of Corollary \ref{indepgens} it is reasonable to  write $\FK{I}$ for $\FK{\AL}$. 
 Since $\PPK{\AL}$ is a projective approximation to $\FK{\AL}$, 
it too depends only on the radical of $I$. We also write 
$K_s^{\bullet}(I) = K_s^{\bullet}(\alpha_1) \otimes \cdots \otimes K_s^{\bullet}(\alpha_n )$, but 
this is an abuse of notation since even its homology groups do depend on the choice of generators. 

\subsection{Local homology and cohomology.}
\label{subsec:lochomcoh}
The local cohomology and homology of an $\AR $-module $M$ are then defined by 
$$ H^*_I(\AR ;M) = H^*(\PPK{I} \otimes M) $$ and
$$   H_*^I(\AR ;M) = H_*(\Hom (\PPK{I} , M). $$
Note that we could equally well use the flat stable Koszul complex in the definition of local 
cohomology, as is more usual. Lemma \ref{Ksupp}
 shows that $H^*_I(M)[1/\beta]=0$ if $\beta \in I$, so $H^*_I(M)$ 
is an $I$-power torsion module and  supported over $I$.

It is immediate from the definitions that local cohomology and local homology are related by 
 a third quadrant universal coefficient spectral sequence 
\begin{equation} 
E_2^{s,t}= \mbox{Ext}^s_\AR (H_I^{-t}(\AR ),M) \Longrightarrow  H_{-t-s}^I(\AR ;M), 
\end{equation} 
with differentials $d_r:E_r^{s,t} \darrow  E_r^{s+r,t-r+1}$. 
%It is observed in \cite{GML} that this specializes to give Grothendieck's 
%local duality spectral sequence \cite{G}.

We observe that local cohomology and homology are invariant under change of base ring. 

\begin{lemma}
If $\AR \lra \AR'$ is a ring homomorphism, $I'$ is the ideal $I\cdot \AR'$ and $M'$ is an 
$\AR'$-module regarded by pullback as an $\AR $-module, then 
$$ \hspace{6mm} H^*_I(\AR ;M')\iso H^*_{I'}(\AR';M') \ 
 \tand \  H_*^I(\AR ;M')\iso H_*^{I'}(\AR';M') $$ \end{lemma}
\begin{proof}
The statement for local cohomology is immediate from the definition, since
$R \tensor_R M'=R'\tensor_{R'}M'$. The statement for local homology is similar, 
for example from the explicit projective model for the stable Koszul complex, or
it can be deduced from the local cohomology statement using the Universal Coefficient Theorem.
\end{proof}

In view of this, we usually  omit the ring $\AR $ from the notation unless is is necessary 
for emphasis. 

It is also useful to know that local homology and cohomology are invariant under completion.

\begin{lemma}
Local homology and cohomology are invariant under the completion $M\darrow M_I^{\sm}$ 
of a finitely generated module $M$. 
\end{lemma}

\begin{proof} From the Universal Coefficient Theorem  it suffices to prove the result for
local cohomology.

First note that the kernel of the natural map $K^{\bullet} (\AL) \lra R$ has
a finite filtration with subquotients $R[1/\beta ]$ for $\beta \in I$. 
Similarly,  the kernel of the natural map $PK^{\bullet} (\AL) \lra R$ has
a finite filtration with subquotients $PR[1/\beta ]$ for $\beta \in I$. 
Since $\beta $ is invertible on $\Hom_R (PR[1/\beta],M)$, Lemma
\ref{Ksupp} shows  $\Hom_R (PR[1/\beta],M) \otimes_R K^{\bullet}(\AL)$ is contractible.

It follows that the natural map  
$$M \otimes_R K^{\bullet}(\AL)
\stackrel{\simeq}\lra   \Hom_R (P K^{\bullet} (\AL), M) \otimes_R K^{\bullet}(\AL)$$
is a homology isomorphism. We will see in \ref{lochomL} below that 
$\Hom_R (P K^{\bullet} (\AL), M) \simeq M_I^{\wedge}$ for finitely generated modules
$M$.
\end{proof}

\subsection{Derived functors}
\label{subsec:derived}

We gave our definitions in terms of specific chain complexes. 
The meaning of the definitions appears in the following two 
theorems.

\begin{thm}[Grothendieck \cite{G}]
\label{loccohR}
If $\AR $ is Noetherian, then the local cohomology
groups calculate the right derived functors of the left exact functor 
$M \longmapsto \Gamma_I (M)$. In symbols, 
$$\hspace{6mm} H_I^n(\AR ;M) = (R^n\Gamma_I)(M). \qqed$$ 
\end{thm}

This result may be used to give an explicit expression for 
local cohomology in familiar terms.  Indeed, 
since $\Gamma_I (M)=\colim_r\, \Hom(\AR /I^r,M)$, and the right
derived functors of the right-hand side are obvious, we
have
$$(R^n\Gamma_I)(M)\iso \colim_r\,\Ext_\AR^n(\AR /I^r,M).$$
The description in terms of the stable Koszul complex is usually 
more practical.

\begin{thm}[Greenlees-May \cite{GML}]
\label{lochomL}
If $\AR $ is Noetherian, then the local homology
groups calculate the left derived functors of the (not usually right exact) 
$I$-adic completion functor $M \longmapsto M_I^{\sm }$. Writing $L_n^I$ for
the left derived functors of $I$-adic completion, this gives 
$$ \hspace{6mm} H_n^I(\AR ;M) = L_n^I(M). \qqed$$ 
\end{thm}

The conclusions of Theorem \ref{loccohR} and  \ref{lochomL} 
are true under much weaker hypotheses \cite{GML,Schenzel}. 

An elementary proof of Theorem \ref{loccohR} can be obtained by induction on 
the number of generators of $I$. To give the idea, we prove the principal case 
$I=(\alpha)$. It is clear that $H^0_I(Q)=\Gamma_I(Q)$ so it suffices
to prove $H^1_I(Q)=\cok(Q \lra Q[1/\alpha])=0$ when $Q$ is an injective $R$-module.

Since $R$ is Noetherian, there is an $r$ for which 
$\alpha^{r+1}x=0$ implies $\alpha^rx=0$, and hence $R \lra (\alpha ) \oplus R/(\alpha^r)$
is a monomorphism. By injectivity of $Q$, any map $q: R \lra Q$ extends to a
 map 
$\langle f,g\rangle: (\alpha) \oplus R/(\alpha^r)\lra Q$, and hence $q=f(\alpha )+g(1)$.
On the other hand, again by injectivity of $Q$, the map $f: (\alpha ) \lra Q$ extends
 to a map $f': R \lra Q$, so $f(\alpha)=\alpha f'(1)$ is a multiple of $\alpha$, and we see
$q=\alpha q'+g(1)$ with $g(1)$ annihilated by $\alpha^r$. This shows 
$$Q=\alpha Q + \Gamma_{\alpha}Q$$
so that $Q \lra Q[1/\alpha]$ is an epimorphism as required. 

The general case  uses the spectral sequence 
$$H^*_I(H^*_J(M))\Longrightarrow H^*_{I+J}(M)$$ 
that is obtained from the isomorphism  
$\PPK{I+J} \cong \PPK{I} \otimes \PPK{J}$. Since we have proved the
 principal case, the remaining step required for the
general case (which we omit) is to show that if $Q$ is injective then $\Gamma_IQ$ is also 
injective. 

The proof of Theorem \ref{lochomL} can also be obtained like this, although it is 
more complicated 
because the completion of a projective module is rarely projective.

\subsection{The shape of local cohomology.}
One is used to the idea that $I$-adic completion is often exact, so that $L^I_0$ is the 
most significant of the left derived functors. However, it is the top non-vanishing 
right derived functor of $\Gamma_I$ that is the 
most significant. Some idea of the shape of these derived functors can be obtained 
from the following result. Observe that the complex  $\PPK{\AL }$ is non-zero only in 
cohomological degrees between $0$ and $n$, so that local homology and 
cohomology are zero above dimension $n$. A result of Grothendieck usually gives 
a much better bound. We write $\dim (R)$ for the Krull dimension of $R$ and
 $\depth_I(M)$ for the $I$-depth of a module $M$
(the length of the longest $M$-regular sequence from $I$).

\begin{thm}[Grothendieck \cite{Go}]
\label{loccohvan}
If $\AR $ is Noetherian of Krull dimension $d$, then 
$$H^i_I(M)=0 \ \  \tand \ \ H_i^I(M)=0 \ \ \text{if} \ i>d.$$ 
If $e=\depth_I(M)$ then 
$$H^i_I(M)=0 \ \text{if} \ i< e.$$
If $\AR $ is Noetherian, $M$ is finitely generated, and $IM\neq M$, then 
$$H^e_I(M)\neq 0. \qqed$$ 
\end{thm}

Grothendieck's proof of  vanishing begins by noting that
 local cohomology is sheaf cohomology with support. It
then proceeds by induction on the Krull dimension and 
reduction to the irreducible case. The statement about depth is 
elementary, and proved  by induction on the length of the $I$-sequence
(see \cite[16.8]{Mat}).

The Universal Coefficient Theorem gives a useful consequence for local 
homology.

\begin{cor}\label{cmcase} 
If $\AR $ is Noetherian and $\mathrm{depth}_I(\AR )=\mathrm{dim}(\AR )=d$, then 
$$\hspace{6mm} L^I_sM=\mathrm{Ext}^{d-s}_\AR (H^d_I(\AR ), M). \qqed $$ 
\end{cor}

For example if $\AR =\Z$ and $I=(p)$, then $H^*_{(p)}(\Z )=H^1_{(p)}(\Z )=\Z /p^{\infty}$. 
Therefore the corollary states that $$L_0^{(p)}M=\mathrm{Ext}(\Z /p^{\infty} , M) \ \ \ \mbox{ and } \ \ \ 
L_1^{(p)}M=\mathrm{Hom}(\Z /p^{\infty} , M), $$
as was observed in Bousfield-Kan \cite[VI.2.1]{BK}.

\subsection{\v{C}ech homology and cohomology.}
\label{subsec:Cech}
We have motivated local cohomology in terms of $I$-power torsion, 
and it is natural to consider the difference between the torsion 
and the original module.  
In geometry this difference would be more fundamental than the torsion 
itself, and local cohomology would then arise by considering functions with support. 

To construct a good model for this difference, 
observe that $\epz: \FK{\AL } \darrow  \AR  $ is an 
isomorphism in degree zero and define the flat \v{C}ech complex $\FC{I}$ to be the complex 
$\SI(\ker\epz)$. Thus, if $i\geq 0$, then $\check{C}^i(I)=K^{i+1}(I)$. For example, if 
$I=(\alpha , \beta )$, then
$$ \FC{I}=\left( \: \AR [1/\alpha ] \oplus \AR [1/\beta ] \darrow  
\AR [1/( \alpha \beta )]\;\right) .$$ 
The differential $K^0(I)\darrow K^1(I)$ specifies a chain map $\AR \darrow \FC{I}$ whose 
fibre is exactly $\FK{I}$. Thus we have a fibre sequence 
$$ \fbox{$\FK{I}  \darrow  \AR  \darrow  \FC{I} .$}$$ 
We define the projective version $\PC{I}$ similarly, using the kernel of the composite of 
$\epz$ and the quasi-isomorphism $\PPK{I}\darrow\FK{I}$; note that $\PC{I}$ is non-zero in 
cohomological degree $-1$.

The \v{C}ech cohomology and homology of an $\AR $-module $M$ are then defined by 
$$ \check{C}H^*_I(\AR ;M) = H^*(\PC{I} \otimes M) $$ 
and 
$$  \check{C}H_*^I(\AR ;M) = H_*(\Hom (\PC{I} , M)). $$ 
The \v{C}ech cohomology can also be defined by use of the flat \v{C}ech complex and is zero 
in negative degrees, but the \v{C}ech homology is usually non-zero in degree $-1$.

The fibre sequence  $\PPK{I} \darrow \AR \darrow \PC{I}$ gives rise to long exact sequences 
relating local and \v{C}ech homology and cohomology, 
$$0 \darrow  H_I^0(M) \darrow  M \darrow  \check{C}H_I^0(M)\darrow  H^1_I(M) \darrow  0$$ and 
$$0 \darrow   H_1^I(M) \darrow   \check{C}H^I_{0}(M)  \darrow  M \darrow  H_0^I(M) 
\darrow   \check{C}H^I_{-1}(M)\darrow  0,$$
together with isomorphisms
$$H^i_I(M) \cong \check{C}H_I^{i-1}(M)
\mbox{ and } H_i^I(M) \cong \check{C}H^I_{i-1}(M)
\mbox{ for } i \geq 2. $$

\subsection{\v{Cech} cohomology and \v{C}ech covers.}
To explain why $\FC{I}$ is called the \v{C}ech complex, we describe how it 
arises by using the \v{C}ech construction to calculate cohomology from a suitable open cover. 
More precisely, let $Y$ be the closed subscheme of $X=Spec(\AR )$ determined by $I$. The space 
$V(I)=\{ \wp | \wp\supset I\}$ decomposes as 
$V(I)=V(\alpha_1 ) \cap \ldots \cap V( \alpha_n )$, and there results an open cover of the 
open subscheme $X-Y$ as the union of the complements 
$X-Y_i$ of the closed subschemes $Y_i$ determined by the principal ideals $(\alpha_i)$. 
However, $X-Y_i$ is isomorphic to the affine scheme $Spec(\AR [1/\alpha_i ])$. Since affine 
schemes have no higher cohomology, 
$$H^*(Spec(\AR [1/\alpha_i ]);\tilde{M})=H^0(Spec(\AR [1/\alpha_i ]);\tilde{M})=M [1/\alpha_i ],$$
where $\tilde{M}$ is the sheaf associated to the $\AR$-module $M$. 
Thus the $E_1$ term of the Mayer-Vietoris spectral sequence for this cover 
collapses to the chain  complex $\FC{I}$, and
$$H^*(X-Y;\tilde{M})\iso \check{C}H^*_I(M).$$

\section{Homotopical analogues of the algebraic definitions} 
\label{sectoploccoh}

We now transpose the algebra from Section \ref{secloccoh} into the homotopy theoretic context.
It is convenient to note that it is routine to extend the algebra to graded rings, 
and we will use this without further comment below. 
We assume the reader is already comfortable working with ring spectra, but there is
an introduction with full references in \cite{spectra}, elsewhere in this 
volume. 

We replace the standing assumption that $\AR$ is a
 commutative $\Z$-algebra  by the assumption that it is 
a commutative $\bbS$-algebra, where $\bbS$ is the sphere 
spectrum. The category of $\AR $-modules is now   the category of 
$R$-module {\em spectra}. Since the derived category of a ring $R$ is equivalent to 
the derived category of the associated Eilenberg-MacLane spectrum \cite{Shipley}, 
the work of earlier sections can be reinterpreted in the new context. 
To emphasize the algebraic analogy, 
 we write $\tensor_R$ and $\Hom_R$ for the smash product over $R$ and function spectrum 
of $R$-maps  and $0$ for the trivial module. In particular $X\tensor_{\bbS}Y=X\sm Y$ 
and $\Hom_{\bbS}(X,Y)=F(X,Y)$.

This section is based on \cite{KEG,Tateca}.

\subsection{Koszul spectra.}
For $\alpha \in \pi_{*} R$, we define the stable 
Koszul spectrum $K(\alpha)$ by the fibre sequence 
$$K(\alpha)\longrightarrow R \longrightarrow R [1/\alpha], $$ 
where 
$R[1/ \alpha] = \hocolim (R \stackrel{\alpha}{\longrightarrow} 
R \stackrel{\alpha}{\longrightarrow} \ldots)$. Analogous to the filtration 
by degree in chain complexes, we obtain a filtration of the 
$R$-module $K(\alpha)$ by viewing it as 
$$\Sigma^{-1}(R[1/\alpha] \cup CR).$$

Next we define the stable Koszul spectrum for the sequence $\alpha_1, \ldots , \alpha_n$ by 
$$K(\alpha_1, \ldots , \alpha_n) = K(\alpha_1) \tensorR \cdots \tensorR K(\alpha_n),$$
and give it the tensor product filtration.

The topological analogue of Lemma \ref{Ksupp} states that if $\beta \in I$ then 
$$K(\alpha_1, \ldots , \alpha_n)[1/\beta ] \simeq 0;$$ 
this follows from Lemma \ref{Ksupp} and the spectral sequence (\ref{algtopss}) below.
We may now use precisely the same proof as in the algebraic case to conclude that the homotopy type of $K(\alpha_1, \ldots , \alpha_n)$ depends only on the radical of the ideal $I = (\alpha_1, \cdots , \alpha_n)$. We therefore write $K(I)$ for 
$K(\alpha_1, \ldots , \alpha_n)$.

\subsection{Localization and completion.}
\label{subsec:loccomp}
With motivation from Theorems \ref{loccohR} and \ref{lochomL}, we define the homotopical 
$I$-power torsion 
(or local cohomology) and homotopical completion (or local homology) modules associated 
to an $R$-module $M$ by
\begin{equation}
\Gamma_I(M)= K(I)\tensorR M \ \ \ \tand \ \ \ \Lambda_I(M)=M^{\wedge}_I=\HomR(K(I), M) . \end{equation} 
In particular, $\Gamma_I(R)=K(I)$.

Because the construction follows the algebra so precisely, it is easy give methods of calculation for the homotopy groups of these $R$-modules. We use the product of the filtrations of the $K(\alpha_i)$ given 
above and obtain spectral sequences
\begin{equation} \label{algtopss}
E^2_{s,t} = H^{-s,-t}_{I} (R_{* }; M_{* }) \Rightarrow \pi_{s+t} (\Gamma_IM) \end{equation} with differentials $d^r: E^r_{s,t} \rightarrow E^r_{s-r,t+r-1}$ and 
\begin{equation}
E_2^{s,t} = H_{-s,-t}^{I} (R^{* }; M^{* }) \Rightarrow \pi_{-(s+t)} (M^{\wedge}_I) \end{equation} with differentials $d_r: E_r^{s,t} \rightarrow E_r^{s+r,t-r+1}.$

\subsection{The \v{C}ech spectra.}
Similarly, we define the \v{C}ech spectrum by the cofibre sequence 
\begin{equation}\label{topKAC}
\fbox{$K(I) \darrow  R \darrow  \check{C}(I). $}
\end{equation}

\noindent We define the homotopical localization (or \v{C}ech 
cohomology) and \v{C}ech homology modules 
associated to an $R$-module $M$ by
\begin{equation}
M[I^{-1}]= \check{C}(I) \tensorR M \ \ \ \mbox{ and }\ \ \ 
 \Delta^I(M)=\HomR(\check{C}(I), M).
\end{equation}
In particular, $R[I^{-1}]= \check{C}(I)$. Once again, we have spectral sequences for 
calculating their homotopy groups from the analogous algebraic constructions.

\subsection{Basic properties.}
We can now give topological analogues of some basic pieces of algebra that we used in 
Section \ref{secloccoh}.  Recall that the algebraic Koszul complex $\FK{I}$ is a direct limit of 
unstable complexes $K_s^{\bullet}(I)$ that are finite complexes of free modules with 
homology annihilated by a power of $I$. We say 
that an $R$-module $M$ is a $I$-power torsion module if its $R_*$-module $M_*$ of homotopy groups 
is a $I$-power torsion module; equivalently, $M_*$ must have support over $I$. 

\begin{lemma} \label{HJcolim}
The $R$-module $K(I)$ is a homotopy direct limit of finite $R$-modules 
$K_s(I)$, each of which has homotopy groups annihilated by some power of $I$. Therefore 
$K(I)$ is a $I$-power torsion module.
\end{lemma}
\begin{proof} It is enough to establish the result for a principal ideal and then 
take tensor products over $R$. Exactly as in the algebraic context, 
since $K(\alpha)$ is the fibre of $R \lra R[1/\alpha]$, 
we find $K(\alpha)\simeq \hocolim_s K_s(\alpha)$, where 
$K_s(\alpha)=\Sigma^{-1} R/\alpha^s$ is the fibre of $\alpha^s :R \darrow  R$.
\end{proof}

The following lemma is an analogue of the fact that $\FC{I}$ is a chain complex which 
is a finite sum of modules $R[1/\alpha ]$ for $\alpha \in I$. 

\begin{lemma} \label{CJfilt}
The $R$-module $\CJI$ has a finite filtration by $R$-submodules with subquotients that are suspensions of modules of the form $R[1/\alpha ]$ 
with $\alpha \in I$. 
\end{lemma}

\begin{proof} By construction $K(\alpha)=\Sigma^{-1}R[1/\alpha]\cup CR$, giving it
 a filtration of length 1. The tensor product $K(I)$ therefore has a filtration 
of length $n$ with the top quotient $R$ and all other subquotients of the form
$R[1/\alpha ]$ for $\alpha \in I$. The result follows since 
 $\CJI$ is the mapping cone of $K(I) \lra R$.
\end{proof}

These lemmas are useful in  combination.

\begin{cor} \label{comp}
If $M$ is a $I$-power torsion module
then $M \tensorR \CJI \simeq 0$; in particular $K(I) \tensorR \CJI \simeq 0$. 
\end{cor} 
\begin{proof} Since $M[1/\alpha ] \simeq 0$ for $\alpha \in I$, Lemma \ref{CJfilt} gives the 
conclusion for $M$. 
\end{proof}

\section{Completion at ideals and Bousfield localization} \label{seccompletionssht}

Bousfield localizations include both completions at ideals and localizations at 
multiplicatively 
closed sets, but one may view these Bousfield localizations as falling into the types 
typified by completion at $p$ and localization away from $p$. Thinking in terms of 
$Spec(R_*)$, 
this is best viewed as the distinction between localization at a closed set 
and localization at the complementary open subset. In this section 
we deal  with the closed sets and  with the open sets in Section \ref{sec:locaway}. 
The section is based on \cite{GM, GML,Tateca}.

\subsection{Homotopical completion.}
\label{htpicalcompletion}
As observed in the proof of Lemma \ref{HJcolim}, we have 
$K(\alpha )=\hocolim_s\Sigma^{-1} R/\alpha^s$ and therefore 
$$M^{\wedge}_{(\alpha)} = \HomR(\hocolim_s \Sigma^{-1} R/\alpha^s , M) 
\simeq \holim_s M/\alpha^s. $$ 
If $I=( \alpha, \beta )$, then 
$$\begin{array}{rcl}
 M^{\wedge}_I &=& \HomR(K(\alpha) \tensorR  K(\beta), M)\\
              &=&\HomR(K(\alpha ), \HomR( K(\beta ), M))\\
              &=&(M_{(\beta) }^{\wedge})_{(\alpha)}^{\wedge},
\end{array}$$
and so on inductively. This should help justify the notation $M_I^{\wedge}=\HomR(K(I),M).$

When $R=\bbS$ is the sphere spectrum and $p \in \Z \cong \pi_0(\bbS)$, $K(p)$ is a Moore 
spectrum for $\Z / p^{\infty}$ in degree $-1$ and we recover the usual definition 
$$X^{\wedge}_p=F(S^{-1}/p^{\infty}, X)$$ 
of $p$-completions of spectra as a special case, where $F(A,B)=\Hom_{\bbS}(A,B)$ is the 
function spectrum. The standard short exact sequence 
for the calculation of the homotopy groups of $X^{\wedge}_p$ in terms of `Ext completion' 
and `Hom completion' follows directly from Corollary \ref{cmcase}.

Since $p$-completion has long been understood to be an example of a Bousfield localization, 
our next task is to show that completion at $I$ is a Bousfield localization in general.

\subsection{Bousfield's terminology.}
 Fix an $R$-module  $E$. A spectrum  $A$ is {\em $E$-acyclic} if $A  \tensorR  E \simeq 0$; 
a map $f:X \darrow  Y$ is an { \em $E$-equivalence} if its cofibre is $E$-acyclic.  
An $R$-module $M$  is {\em $E$-local} if $E \tensorR T \simeq 0$ 
implies  $\HomR (T, M) \simeq 0$. A map $Y \darrow  L_EY$ is a 
{\em Bousfield $E$-localization} of  $Y$ if it is an $E$-equivalence 
and $L_EY$   is  $E$-local. This means that $Y \darrow  L_EY$ is  terminal
among $E$-equivalences with domain $Y$, and the Bousfield localization is therefore unique 
if it exists. Similarly, we may replace the single spectrum $E$ by a class $\sE$ 
of objects $E$, and require the conditions hold for all such $E$ 

The following is a specialization of a change  of rings result to the
ring map $\bbS \lra R$.

\begin{lemma} \label{loclem}
Let ${\sE}$ be a class of $R$-modules. If an $R$-module $N$ is 
${\sE}$-local as an $R$-module, then it is ${\sE}$-local as an $\bbS$-module. 
\end{lemma} 
\begin{proof} If $E\sm T =E\tensor_{\bbS} T \htp *$ for all $E$, 
then $E\tensorR (R \tensor_{\bbS} T)\htp 0$ for all $E$ and therefore 
$F(T,N)=\Hom_{\bbS}(T,N) \htp \HomR(R \tensor_{\bbS} T,N) \htp 0$. 
\end{proof}

\subsection{Homotopical completion is a Bousfield localization.}
The  class that will concern us most is the class $I$-{\bf Tors} of finite $I$-power torsion $R$-modules $M$. Thus $M$ must be a finite cell $R$-module, and its $R_*$-module $M_*$ of homotopy groups must be a $I$-power 
torsion module. 

\begin{thm} \label{Bousclosed}
For any finitely generated ideal $I$ of $R_*$ the map $M \darrow  M_I^{\wedge}$ is 
Bousfield localization in the category of $R$-modules in each of the following 
equivalent senses:
\begin{enumerate}
\item[(i)] with respect to the $R$-module $\Gamma_I(R)=K(I)$. \item[(ii)] with respect to the class $I$-{\bf Tors} of finite $I$-power torsion $R$-modules. 
\item[(iii)] with respect to the $R$-module $K_s(I)$ for any $s \geq 1$. 
\end{enumerate} Furthermore, the homotopy groups of the completion are related to local homology groups 
by a spectral sequence 
$$E^2_{s,t}=H^I_{s,t}(M_*)\Longrightarrow \pi_{s+t}(M_I^{\wedge}).$$ 
If $R_*$ is Noetherian, the $E^2$ term consists of the left derived functors of $I$-adic completion: $H^I_s(M_*)=L_s^I(M_*)$. \end{thm} 
\begin{proof}  We begin with (i). Since 
$$\HomR(T, M_I^{\wedge}) \simeq \HomR(T \tensorR K(I), M),$$
it is immediate that $M_I^{\wedge}$  is $K(I)$-local. We must prove that the map $M \darrow  M_I^{\wedge}$ is a $K(I)$-equivalence. The fibre of this map is $\Hom_R (\CJI , M)$, so we must show that  
$$\Hom_R (\CJI , M) \tensorR K(I) \simeq 0.$$
By Lemma \ref{HJcolim}, $K(I)$ is a homotopy direct limit of terms $K_s(I)$. Each $K_s(I)$ is in $I$-{\bf Tors}, and we see by their definition in terms of cofibre sequences and smash products that their duals 
$K_s(I)^{\#}$ are also in $I$-{\bf Tors}, where $M^{\#}=\HomR(M,R)$. 
Since $K_s(I)$ is a finite cell $R$-module,
$$\HomR(\CJI , M) \tensor_R K_s(I)= \HomR( \CJI \tensorR K_s(I)^{\#},M),$$ and $\CJI \tensorR K_s(I)^{\#} \simeq 0 $ by Corollary \ref{comp}. Parts (ii) and 
(iii) are similar but simpler. For (iii), observe that we have a cofibre sequence 
$R/\alpha^s \darrow  R/\alpha^{2s} \darrow  R/\alpha^s$, so that all of the $K_{js}(I)$ may be constructed from $K_s(I)$ using a finite number of cofibre sequences. 
\end{proof}

\section{Localization away from ideals and Bousfield localization}
\label{sec:locaway}

In this section we turn to localization away from the closed set defined
by an ideal $I$. First, observe that, when $I=(\alpha)$, 
$M[I^{-1}]$ is just $R[\alpha^{-1}]\tensor_R M = M[\alpha^{-1}]$. However, the higher 
\v{C}ech cohomology groups give the construction for general finitely generated ideals 
a quite different algebraic flavour, and $M[I^{-1}]$ is rarely a localization 
of $M$ at a multiplicatively closed subset of $R_*$. 
The section is based on \cite{Tateca}.

\subsection{The \v{C}ech complex as a Bousfield localization.}
To characterize this construction as a Bousfield 
localization, we consider the class $I$-{\bf Inv} of $R$-modules $M$ for which there is 
an element $\alpha \in I$ such that $\alpha :M \darrow M$ is an equivalence. 

\begin{thm}
For any finitely generated ideal $I=(\alpha_1,\ldots,\alpha_n)$ of $R_*$, the map $M \darrow M[I^{-1}]$ is Bousfield localization in the category of $R$-modules in each of the following equivalent senses: \begin{enumerate} \item[(i)] with respect to the $R$-module $R[I^{-1}]=\CJI$. \item[(ii)] with respect to the class  $I$-{\bf Inv}. \item[(iii)] with respect to the set $\{ R[1/\alpha_1 ], \ldots , R[1/\alpha_n] \}$. \end{enumerate} Furthermore, the homotopy groups of the localization are related to \v{C}ech cohomology 
groups by a spectral sequence
$$E_{s,t}^2 = \check{C}H_I^{-s,-t}(M_*) \Longrightarrow \pi_{s+t}(\CIM ).$$ If $R_*$ is Noetherian, the $E^2$ term can be viewed as the cohomology of
$Spec(R_*) \setminus V(I)$ with coefficients in the sheaf associated to $M_*$. 
\end{thm} 
\begin{proof} We prove Part (i); Parts  (ii) and (iii) are proved similarly.
To see that $\CIM$ is local, suppose that $T \tensorR \CJI  \simeq 0$. We must show that 
$\HomR(T, \CIM ) \simeq 0$. By the cofibre sequence defining $\CJI$ and the supposition, it 
suffices to show that 
$\HomR(K(I) \tensorR T, \CIM ) \simeq 0$. By Lemma \ref{HJcolim}, 
$$\HomR(K(I) \tensorR T, \CIM) \simeq \holim_s \HomR(K_s (I) \tensorR T,\CJI \tensorR M).$$ 
Observing that
$$\HomR(K_s(I) \tensorR T, \CJI  \tensorR M) \simeq \HomR( T, K_s(I)^{\#} \tensorR \CJI \tensorR M),$$ 
we see that the conclusion follows from Corollary \ref{comp}. The map $M \darrow  \CIM$ is a 
$\CJI$-equivalence since its fibre is $\Gamma_I(M)=K(I) \tensorR M$ and 
$K(I)\tensorR \CJI  \simeq 0$ by Corollary \ref{comp}.
 \end{proof}

\begin{remark}
Translating the usual terminology, we say that a localization $L$ on $R$-modules 
is {\em  smashing} if $L(N)=N \tensorR L(R)$ for all $R$-modules  $N$. 
It is clear that localization away from $I$ is smashing and  that completion 
at $I$ will usually not be smashing. For any smashing localization
one may form an associated complementary completion, giving a formal situation
like the present one \cite{axiomatic}.
\end{remark}

\begin{remark}
\label{rk:cellular}
One may also characterize the map $\Gamma_I (M) \lra M$ by a universal property
analogous to that of the cellular approximation in spaces. 

On the one hand, $\Gamma_I(M)$ is constructed from $K_1(I)$ by \ref{HJcolim},  
and on the other hand, the map induces an equivalence of $\HomR (K_1(I), \cdot )$
since, by Lemma \ref{CJfilt}, 
$\HomR(K_1(I), \CJI )\simeq 0.$
We will return to these ideas in Part 2.
\end{remark}

\section{Chromatic filtrations for $MU$-modules.}
\label{sec:bordism}

We have described ways of picking out parts of stable homotopy theory
concentrated on open and closed subschemes of a coefficient ring. This 
philosophy was used before in stable homotopy theory in a deeper way,  where 
the analogy needed much more work to set up. We show in this section that 
the two approaches agree where they overlap. 

The starting point of the chromatic apprach to homotopy theory is that 
one may attempt to make calculations in stable homotopy theory using a
chosen cohomology theory. One of the most effective is complex cobordism $MU^*(X)$.
Not only is it a very practical means of calculation in many cases, but 
it also leads to remarkable structural insight \cite{Hopkinssurvey, Rave, Rave2, 
HoveyStrickland}. Indeed, the Adams-Novikov spectral sequence gives a means of calculation 
in  the category  of $\bbS$-modules based on the 
 category of $MU_*$-modules. However the connection goes further. 
The coefficient ring $MU_*$  is isomorphic Lazard's universal ring for 1-dimensional commutative 
formal group laws, and Quillen proved that there is a canonical isomorphism for geometric
reasons. Thus the algebra of $MU_*$-modules is that of 1-dimensional commutative formal 
 groups, and this gives rise to a filtration by height. This algebraic filtration in 
 turn gives rise to the chromatic filtration of stable homotopy.

 Since cobordism is represented by a commutative $\bbS$-algebra  $MU$, 
 there is a well-behaved category of $MU$-modules, and we may lift the Adams-Novikov 
spectral sequence to the 
 level of $MU$-modules. The  formalities we have described in Section \ref{sectoploccoh} 
above give a 
 filtration of the category of $MU$-modules, and it is the purpose of this section 
 to explain that this agrees with the restriction of the chromatic filtration, and
 to provide a dictionary for relating the languages.  

\subsection{Cobordism notation.}
Recall that $MU_*=\Z [x_i\,|\,i\geq 1]$, where $\mathrm{deg} (x_i)=2i$, and that 
 $MU_*$ contains elements 
$v_i$ of degree $2(p^i-1)$ that map to the Hazewinkel generators of 
$BP_*=\Z_{(p)}[v_i\,|\,i\geq 1]$. We let $I_n$ denote the ideal $(v_0,v_1,\ldots ,v_{n-1})$ in 
$\pi_*(MU)$, where $v_0=p$; we work with $MU$ rather than $BP$ because of its canonical 
$\bbS$-algebra structure. The Brown-Peterson espectrum $BP$ is a ring 
in the homotopy category of $MU$-modules whose unit 
$MU\darrow BP$ factors through the canonical retraction $MU_{(p)}\darrow BP$. There are
also $MU$-module spectra $E(n)$  and $K(n)$ such that 
$E(0)_*=K(0)_*=\Q$,  and for $n\geq 1$ 
$$K(n)_*=\Fp [v_n,v_n^{-1}] 
\mbox{ and } E(n)_*=\Bbb{Z}_{(p)}[v_{1},\dots,v_{n},v_{n}^{-1}].$$ 
Furthermore $E(n)$ is a ring  in the homotopy category of $MU$-modules.

\subsection{Chromatic filtration and the Cousin complex.}

The chromatic filtration corresponds to the Bousfield localizations
$L_n$ with respect to $K(0)\vee K(1) \vee \ldots \vee K(n)$, or equivalently 
with respect to $E(n)$. 
For any spectrum $X$, Ravenel defines the $n$th acyclization $C_nX$ and the 
$n$th monochromatic part $M_nX$ by the cofibrations
$$C_nX \lra X \lra L_nX$$
and
$$\Sigma^{-n}M_nX \lra L_n X \lra L_{n-1}X.$$
To  start inductions, $L_{-1}X=*$.  

The connection between the present constructions and Ravenel's
are provided by the following dictionary.

\begin{lemma}
For an $MU$-module $A$, and for any $n \geq -1$,  we have
$$L_{n}A \htp A[I_{n+1}^{-1}],$$ 
$$C_{n}A \htp \Gamma_{I_{n+1}}A$$
and  
$$M_nA \htp \SI^{n-1}\Gamma_{I_n}A[1/v_n].$$ 
\end{lemma}

\begin{proof}
It suffices to establish the first identification. % between the $C$s and $\Gamma$s. 
Indeed,  $C_{n}A$ and $M_{n}A$ are then defined 
by cofibre sequences from  the $C_i$, and these are precise counterparts of 
$$\Gamma_{I_{n+1}}A \lra A \lra A[I_{n+1}^{-1}]$$
and
$$\Sigma^{-1}\Gamma_{I_{n}}A [1/v_n] \lra \Gamma_{I_{n+1}}A \lra \Gamma_{I_{n}}A $$
from Subsection \ref{subsec:Cech}.

The idea of the comparison is clear, but the obstacle is that in considering
$MU_*(X)$ for an $MU$-module $X$ we must relate the $MU_*$-action
through $MU$ to the $MU_*$-action through $X$. We need to engineer a situation where
 we can apply the well-known formulae this purpose. 
By \cite[7.3.2]{Rave2}, localization at $E(n)$ is the same as localization at 
the wedge of $MU[1/v_i]$ for $0 \leq i \leq n$. By \ref{loclem} we conclude that 
$M[I_{n+1}^{-1}]$  is $E(n)$-local. To see that $M \lra M[I_{n+1}^{-1}]$ is an 
$MU[1/v_i]$-equivalence for $0 \leq i \leq n$, note that its fibre is 
$\Gamma_{I_{n+1}}M$ and $\Gamma_{I_{n+1}}M[1/w]\simeq 0$ for any $w \in I_{n+1}$. 
Consider $\pi_*(MU \sm MU)$ as a left $MU_*$-module, and recall from \cite[B.5.15]{Rave2}
that the right unit $\eta_R: MU_* \lra \pi_*(MU \sm MU)$ satisfies
$$\eta_R(v_i)\equiv v_i \mbox{ mod } I_i \cdot \pi_*(MU \sm MU) \mbox{, hence }
\eta_R(v_i) \in I_{i+1}\cdot \pi_*(MU \sm MU).$$
We have
$$\Gamma_{I_{n+1}}(M) \sm MU \simeq \Gamma_{I_{n+1}}(M) \tensor_{MU} (MU \sm MU)$$
and can deduce inductively that $\Gamma_{I_{n+1}}M \sm MU[1/w] \simeq 0$ for 
any $w \in I_{n+1}$ since $\Gamma_{I_{n+1}}(M)[1/w] \simeq 0$ for any such $w$.
\end{proof}

\subsection{Chromatic completions of $MU$-modules and $\bbS$-modules.}

We first recall the definition of the chromatic completion.
For a sequence $\bold{i} = (i_0, i_1, \ldots , i_{n-1})$, 
we may attempt to construct generalized Toda-Smith spectra 
$$M_{\bold{i}}=M(p^{i_0}, v_1^{i_1}, \ldots , v_{n-1}^{i_{n-1}})$$ 
inductively, starting with $\bbS$, continuing with the cofibre sequence  
$$M(p^{i_0}) \darrow  \bbS \stackrel{p^{i_0}} \darrow  \bbS,$$
and, given $L=M_{(i_0, i_1, \ldots , i_{n-2})}$, concluding with the cofibre sequence 
$$M_{\bold{i}} \darrow  L \overto{v_{n-1}^{i_{n-1}}} L.$$ 
Here $M_{\bold{i}}$ is a finite complex of type $n$  and hence admits a $v_n$-self 
map by the Nilpotence Theorem \cite{DHS, HS}, and $v_n^{i_n}$ is shorthand for such a map. 
These spectra do not exist for all sequences {\bf i}, but they do exist for a cofinal set of 
sequences, and Devinatz has shown \cite{D} that there is a cofinal collection all of which 
are 
ring spectra. These spectra are not determined by the sequence, but it follows from the 
Nilpotence Theorem that they are asymptotically unique in the sense that 
$\hocolim_{\bold{i}}  M_{\bold{i}}$ is independent of all choices. Hence we may define a 
completion for all $p$-local spectra $X$ by 
$$X_{I_n}^{\wedge} = F( \hocolim_{\bold{i}}  M_{\bold{i}}, X).$$ 
We denote the spectrum $\hocolim_{\bold{i}}M_{\bold{i}}$ by  $\Gamma_{I_n}(\bbS)$, 
although it is not simply a  local cohomology spectrum. 

\begin{prop} Localize all spectra at $p$. Then there is an equivalence of $MU$-modules 
$$ MU\sm \Gamma_{I_n}(\bbS) \simeq \Gamma_{I_n}(MU).$$
Therefore, for any $MU$-module $M$, there is an equivalence of $MU$-modules between the 
two completions $M_{I_n}^{\wedge}$. 
\end{prop} 
\begin{proof} The second statement follows from 
the first since
$$\Hom_{MU}(MU\sm \Gamma_{I_n}(\bbS),M)\htp F(\Gamma_{I_n}(\bbS),M)$$
as $MU$-modules.
It suffices to construct compatible equivalences
$$MU\sm M_{\bold{i}}\htp MU/p^{i_0} \tensor_{MU} MU/v_1^{i_1} \tensor_{MU}  
\ldots \tensor_{MU} MU/ v_{n-1}^{i_{n-1}}.$$
The right side is equivalent to $MU/I_{\bold{i}}$, where 
$I_{\bold{i}}= (p^{i_0},v_1^{i_1}, \ldots , v_{n-1}^{i_{n-1}})\subset I_n$. A $v_n$-self 
map $v: X\darrow X$ on a type $n$ finite complex $X$ can be characterized as a map such that, for some $i$, $BP_*(v^i): BP_*(X)\darrow BP_*(X)$ is multiplication by $v_n^j$ for some $j$. 
Since $MU_*(X)=MU_*\tensor_{BP_*}MU_*(X)$, we can use $MU$ instead of $BP$. Using $MU$, we 
conclude that the two maps of spectra $\id\sm v^i$ and $v_n^j\sm \id$ from $MU\sm X$ to 
itself induce the same map on homotopy groups. The cofibre of the first is $MU\sm Cv^i$ 
and the cofibre of the second is $MU/(v_n^j)\sm X$. If $X$ is a generalized Moore spectrum, 
results of  \cite{HS} show that some powers of these two maps are homotopic, 
hence the cofibres of these powers are equivalent. The conclusion follows by induction. 
\end{proof}

\section{Completion theorems and their duals.}
\label{sec:compthms1}

Finally, we return to the motivating context with a section 
based on \cite{GM,KEG}. There is more background on equivariant
stable homotopy theory in \cite{Handbook1}. So as to illustrate
the result in the simplest context we choose a finite
group $G$, and an equivariant cohomology theory $R_G^*(\cdot)$
which is Noetherian in the sense that $R_G^*$ is Noetherian and 
for any finite $G$-CW-complex $X$, the $R_G^*$-module $R_G^*(X)$ is 
finitely generated. We also suppose that $R_G^*(X)\cong R^*(X/G)$ if $X$ is $G$-free
(in fact we make the corresponding assumption about the representing $G$-spectrum, which is
to say that  it is split \cite{GMt}). It smooths the exposition to suppose throughout
that $G$-spaces are equipped with a $G$-fixed basepoint. This is no loss of generality 
since an unbased $G$-space $Y$ may be converted to a based $G$-space $Y_+$ by adding a 
disjoint basepoint. The unreduced cohomology of $Y$ is the reduced cohomology of $Y_+$.

The completion theorem concerns the map 
$$R_G^*=R_G^*(S^0) \lra R_G^*(EG_+)=R^*(BG_+)$$
induced by projection $EG \lra pt$, or more generally
$$R_G^*(X) \lra R_G^*(EG_+ \sm X)=R^*(EG_+\sm_G X)$$
for a finite based $G$-CW-complex $X$. 
The classical completion theorem is said to hold if these maps are completion 
at the augmentation ideal 
$$I=\ker(R_G^* \lra R^*).$$

By the Noetherian condition, $I$ is finitely generated, so if 
$R_G^*(\cdot )$ is represented by a strictly commutative ring spectrum, we may form 
the $I$-adic completion of
$R$. The classical completion theorem for the case $X=pt_+=S^0$  
then suggests the derived completion theorem:
$$F(EG_+,R) \simeq R_I^{\wedge}. $$
More precisely, the statement is that the natural map $R \lra F(EG_+,R)$
has the universal property of $I$-adic completion.

If the derived completion theorem holds, the classical completion theorem 
holds by taking homotopy groups, since the higher derived functors vanish 
by the Artin-Rees Lemma. Furthermore, applying $F(X, \cdot)$ 
we obtain 
$$F(X \sm EG_+,R) \simeq F(X,R)_I^{\wedge}. $$
The classical completion theorem for a general finite $X$ follows by 
taking homotopy groups.  Furthermore, we can see how to 
formulate the appropriate statement for an infinite $X$ using local 
homology \cite{GM}

Now notice that $R_I^{\wedge}=\Hom_R(\Gamma_IR,R)$ and 
$F(EG_+,R)\simeq \Hom_R(R \sm EG_+, R)$. Thus  the derived completion 
theorem is the statement that  $\Gamma_IR $ and $R \sm EG_+$ become
equivalent after applying $\Hom_R(\cdot  , R)$. Indeed, there are obvious
comparison maps 
$$R \sm EG_+ \stackrel{\simeq}\lla \Gamma_I(R \sm EG_+) \lra \Gamma_IR, $$
and it is elementary that the left hand map is an equivalence:
 the map $\Gamma_I(X \sm F ) \lra X \sm F$ is an equivalence
when $F=G_+$ since $\res^G_1(I)=(0)$, and hence for $F=EG_+$, since
it is constructed from $G$-free cells. 

The simplest way that the derived completion could be true would be
if the derived local cohomology theorem holds:
$$R \sm EG_+ \simeq \Gamma_IR, $$
or more precisely that $R \sm EG_+ \lra R$ has the universal property 
of $\Gamma_IR \lra R$. 

The remarkable fact is that for complex oriented cohomology theories
this strengthening of the completion theorem holds. Indeed, the 
proof is simple enough that  we can outline it here. However there are examples (such 
as stable cohomotopy) for which the derived completion theorem holds
(by Carlsson's proof \cite{Carlsson} of the Segal Conjecture) but for which the local 
 cohomology theorem is false.

\begin{thm}
For a complex oriented, Noetherian theory, the derived local cohomology 
theorem holds in the sense that
$$R \sm EG_+ \simeq \Gamma_IR. $$
Hence in particular there is a spectral sequence
$$H^*_I(R_G^*)\Rightarrow R^G_*(EG_+)\cong R_*(BG_+). $$
\end{thm}

\begin{proof} 
We need only show that the map $\Gamma_I(R \sm EG_+) \lra \Gamma_I(R)$
is an equivalence, or equivalently that $\Gamma_I(R \sm \tilde{E}G)\simeq 0$, 
where $\tilde{E}G$ is the mapping cone of $EG_+ \lra S^0$.
Since the theory is complex oriented, we may define
Euler clases of complex representations.

If $G$ is of prime order we may take $\tilde{E}G =S^{\infty V}$
where $V$ is the complex reduced regular representation and then 
use the fact that the Euler class $\chi (V) $ is in $I$, combined with 
Lemma \ref{Ksupp}. 

The same argument shows that $\Gamma_I(R \sm S^{\infty V})\simeq 0$
for any complex representtion $V$. We use an induction on the group order, noting that if 
$V$ is the reduced regular representation, there is
 a map $\tilde{E}G \lra S^{\infty V}$ whose mapping cone $C$ is built
from cells $G/H_+$ with $H$ a proper subgroup. 
The Noetherian hypothesis is enough to ensure the radical of 
 $\res^G_H(I(G))$ is $I(H)$, so it follows by induction that
$\Gamma_I(R \sm \tilde{E}G \sm G/H_+) \simeq 0$, and hence
that $\Gamma_I(R \sm \tilde{E}G \sm C) \simeq 0$ as required. 
\end{proof}

\part{Morita equivalences and Gorenstein rings.}
%\section{Introduction to Part 2.}

In Part 1 we considered certain localizations and completions associated
to an ideal $I$ of the coefficient ring $R_*$ of a ring spectrum $R$, 
and gave characterizations up to homotopy, but the methods were based
on using particular elements of $R_*$. In Part 2, we change focus 
from the pair $(I,R)$  to the pair $(R, R_*/I)$, and then replace the graded ring
$R_*/I$ by a ring spectrum $k$. This is a strict refinement
when the coefficient ring of $k$ is the quotient considered before,
in the sense that $k_*=R_*/I$. Conceptually, the present situation 
is considerably more general, 
but for a particular ring spectrum  $R$ and a particular ideal $I$ 
of its coefficients it is sometimes impossible to construct an $R$-module $k$ 
(let alone an $R$-algebra) with $k_*=R_*/I$.

To summarize: in Part 2 we consider a map $R \lra k$ of ring spectra, and
 the notation is chosen to suggest the analogy  with a commutative Noetherian 
local ring $R$ with residue field $k$. In effect we are extending the idea of trying 
to do commutative algebra entirely in the derived category. 
Given the map $R\lra k$, it turns out that completion and localization 
arise from a Morita  adjunction, so that it is essential that we consider
non-commutative rings, even if we only want conclusions in commutative algebra. 
Part 2 is based on the papers \cite{tec} and \cite{DGI}.

\section{The context, and some examples.}
\label{sec:context}

The basic ingredients are as follows.
\begin{context}
The main input is a map $R\lra k$ of ring spectra with notation suggested 
by the case when $R$ is commutative local ring with residue field $k$. 
We also write $\cE =\HomR (k,k)$ for the (derived) endomorphism ring spectrum.
 \end{context}

We work in the derived category $D(\Rmod)$ of left $R$-modules. 

\subsection{New modules from old.}
Three construction principles will be important to us. There
is some duplication in terminology, but the flexibility is convenient.

If $M$ is an $R$-module we say that $X$ is {\em built} from $M$ if $X$ can be 
formed from $M$ by completing triangles, taking 
coproducts and retracts (i.e., $X$ is in the {\em localizing subcategory} 
generated by $M$). We refer to objects built from $M$ as {\em $M$-cellular}, 
and write  $\cell (R,M)$ for the resulting full subcategory of $D(\Rmod )$.
An {\em $M$-cellular approximation} of $X$ is a map $\cell_M(X) \lra X$ where
$\cell_M(X)$ is $M$-cellular and the map is an $\HomR (M,\cdot)$-equivalence.

 We say that it is {\em finitely built} from $M$ if only finitely many steps 
and finite coproducts are necessary (i.e., $X$ is in  the {\em thick} 
subcategory generated by $M$).

Finally, we say that $X$ is {\em cobuilt} from $M$ if $X$ can be 
formed from $M$ by completing triangles, taking 
products and retracts (i.e., $X$ is in the {\em colocalizing subcategory} 
generated by $M$).

\subsection{Finiteness conditions.}

We say $M$ is {\em small} if it is finitely built from $R$
and we say $R$ is {\em regular} if $k$ is small as an $R$-module. 
The theory is easiest if $R$ is regular, but this  assumption is 
very strong. We will develop a useful theory under a much weaker condition. 

\begin{defn} \cite{DGI} 
\label{assumptionKk}
We say that $k$ is {\em proxy-small} if there is an object 
$K$ with the following properties
\bi
\item $K$ is small 
\item  $K$ is finitely built from $k$ and
\item  $k$ is built from $K$.
\ei
\end{defn}

\begin{remark}
Note that the second and third condition imply that the $R$-module $K$ 
generates the same category as $k$ using
triangles and coproducts: $\cell (R,K)=\cell (R,k) . $
\end{remark}

It is one of the messages of \cite{DGI} that the weak condition 
of  proxy-smallness allows one to develop a very useful theory.

\subsection{The principal examples.} We will illustrate the theory
throughout with three rather extensive families of examples, 
giving them in this order each time.

\begin{examples} 

(i) {\em (Local algebra)}
Take $R$ to be a commutative Noetherian local ring in degree 0, with maximal 
ideal $I$ and residue field $k$. 

By the Auslander-Buchsbaum-Serre theorem, $k$ is  small if and only if $R$ is a regular 
local ring,  confirming that  the smallness of $k$ is a very 
strong condition. On the other hand, $k$ is  always proxy-small:  
we may take $K=K_1(\AL)$ to be the Koszul complex for a generating 
sequence $\AL$ for $I$.

It is shown in \cite{tec} that $\cell (R,k)$ consists of objects whose 
homology is $I$-power torsion. 

The endomorphism ring $\cE$ is modelled by the endomorphism ring of 
a projective $R$-resolution of the $R$-module $k$. Its homology is the 
Ext algebra
$$H_*(\cE)=\ExtR^* (k,k).$$

(ii) {\em (Cochains on a space)}
For a field $k$, consider the cochain complex for a space $X$ 
$$R=C^*(X; k). $$ 
We want a commutative model, so it is not possible to use a conventional
differential graded algebra unless $k$ is of characteristic 0. Instead, 
we identify the ring $k$ with the associated Eilenberg-MacLane spectrum 
$Hk$, chosen to be a strictly commutative ring spectrum. 
We will therefore interpret $C^*(X;k)=F(X,Hk)$ as the function spectrum of 
maps from $X$ to the 
commutative ring spectrum $Hk$, and as such it is itself a commutative $k$-algebra
spectrum. This is reasonable since $C^*(X;k)$ is a model for ordinary 
cohomology in the sense that $\pi_*(C^*(X;k))=H^*(X;k)$.

If the Eilenberg-Moore spectral sequence converges then 
$$\cE \simeq  C_*(\Omega X; k). $$
For instance this holds if $X$ is simply connected, or if $k$ 
is characteristic $p$ and  $X$ is connected with 
fundamental group a finite $p$-group \cite{DwyerEM}. 

Under this assumption one may show that $R$ is regular if and only 
if $H_*(\Omega X;k)$ is finite dimensional. In one direction, we
 use the Eilenberg-Moore construction:  if $k$ is finitely built from 
$R$ then we apply $\HomR (\cdot , k)$ to deduce
$\HomR (k,k)\simeq C_*(\Omega X;k)$ is finitely built from 
$\HomR (R,k)=k$. The other direction uses the Rothenberg-Steenrod 
construction.  The ring $R$ is proxy-regular much more generally. 
In particular $R$ is proxy-regular if $X=BG$ for a compact Lie group 
$G$ or if $X$ is a connected finite complex, whatever the fundamental 
groups.

(iii) {\em (Equivariant cohomology)}
To make contact with the equivariant motivation in Part 1
we suppose given a commutative ring $G$-spectrum, $R$ representing
$R_G^*(X )=[X,R]_G^*$, and consider $k=F(G_+,R)$ representing the non-equivariant 
theory $R^*(X)=[X,R]^*=[X,F(G_+,R)]_G^*$. The collapse
map $G_+ \lra S^0$ induces a ring map $R \lra k$ and 
$$\cE = \HomR (k,k)\simeq \HomR (F(G_+,R),F(G_+,R)) \simeq F(G_+,G_+) \sm R.$$
The category $\cell (R,k)$ consists of  the $G$-free $R$-modules.
\end{examples}

\section{Morita equivalences.}
\label{sec:Morita}

Morita theory studies objects $X$ of a
category $\C$ by considering their counterparts $\Hom (k,X)$
as modules over the endomorphism ring $\End (k)$ of an object $k$.
In favourable circumstances this may provide an equivalence
between $\C$ and a category of $\End (k)$-modules. 
Classically,  $\C$ is an abelian category with infinite sums and $k$ 
is a small projective generator, and we find $\C$ is equivalent
to the category of $\End(k)$-modules \cite[II Thm 1.3]{BassK}.  However there are numerous
variations and extensions. The differences in our context that are
 most significant are that we work in a derived category (or its model) 
rather than in an 
abelian category, and that $k$ is not necessarily either small or
a generator. The fact that the objects of the categories are
spectra is unimportant except for the formal context it provides.
See \cite{SchwedeM} for an account from the present point of view.

 Two separate Morita equivalences play a role. It is striking that
even where we want conclusions in a commutative context, the
Morita equivalences involve us in thinking about non-commutative
rings. This may seem less surprising if we think of it as a 
refinement of  the well established use of non-commutative rings 
of operations in  stable homotopy theory.
This section is based on \cite{tec}, with augmentations from \cite{DGI}.

\subsection{First variant.}
Continuing to let $\cE=\HombR (\boldk,\boldk)$ denote the (derived) 
endomorphism ring, we consider the relationship between the derived 
categories of left $\boldR$-modules 
and  of right $\cE$-modules. We have  the adjoint pair 
$$\adjunction{T}{D(\modcE)}{D(\bRmod)}{E}$$
defined by 
$$T(X):=X \tensorcE \boldk \mbox{ and } E(M):=\HombR (\boldk,M).$$ 

\begin{remark}
 If $k$ is small, it is easy to see that this adjunction gives equivalence 
$$\cell (R,k) \simeq D(\modcE)$$
between the derived 
category of $\boldR$-modules built from $\boldk$
and the derived category of $\cE$-modules. Indeed, to see the 
unit $X \lra ETX=\HomR (k, X \tensorE k)$ is an equivalence, 
we note it is obviously an equivalence  for $X=\cE$ and hence
for any $X$ built from $\cE$, by smallness of $k$. The argument
for the counit is similar.
\end{remark}

\begin{remark}
The unit of the adjunction is not an equivalence in general. For 
example if $R=\Lambda (\tau)$ is exterior on a generator of degree $1$ then 
$\cE \simeq k[x]$ is polynomial on a generator of 
degree $-2$. As an $R$-module,  $k$ is of infinite homological dimension and
hence it is not small. In this case all $R$-modules are $k$-cellular, so that 
$\cell (R,k) =\Rmod$. Furthermore, the only subcategories of $R$-modules 
closed under coproducts
and triangles are the trivial category and the whole category. 
On the other hand the category of torsion $\cE$-modules is a proper 
non-trivial subcategory closed under coproducts and triangles. 

Exchanging roles of the rings, so that $R=k[x]$ and $\cE \simeq \Lambda (\tau)$, 
we see $k$ is small as a $k[x]$-module 
and $\cell (k[x],k)$ consists of torsion modules. Thus we deduce
$$\mbox{tors-$k[x]$-mod} \simeq \mbox{mod-$\Lambda (\tau)$}. \qqed$$
\end{remark}

In fact  the counit 
$$TEM =\HomR (k,M) \tensorE k \lra M$$ 
of the adjunction is of interest much more generally. Notice that 
any $\cE$-module (such as $\HomR (k,M)$) is built from 
$\cE$, so the domain is $k$-cellular. We say $M$ is {\em effectively 
constructible} from $k$ if the counit is an equivalence, because
$TEM$ gives a concrete and functorial model for the cellular approximation
to $M$.  Under the much weaker assumption of proxy smallness 
we obtain a very useful conclusion linking Morita theory to commutative
algebra. 

\begin{lemma}
\label{proxycell}  
Provided $\boldk$ is proxy-small, 
 the counit 
$$TEM =\HomR (k,M) \tensorE k \lra M$$ 
is $k$-cellular approximation, and hence
in particular any $k$-cellular object is effectively constructible from $k$. 
\end{lemma}

\begin{proof} 
We observed above that the domain 
is $k$-cellular. To see the counit is a $\HomR (k,\cdot )$-equivalence, 
consider the evaluation map 
$$\gamma: \HomR (k,X) \tensorE \HomR (Y,k) \lra \HomR (Y,X).$$ 
This is an equivalence if $Y=k$, and hence by proxy-smallness it is an 
equivalence if $Y=K$. This shows that the top horizontal in the  diagram
$$\begin{array}{ccc}
\HomR (k,X) \tensorE \HomR (K,k)& \stackrel{\simeq} \lra &\HomR (K,X)\\
\simeq \downarrow &&\downarrow = \\
\HomR (K,\HomR (k,X)\tensorE k) &\lra &\HomR (K,X)
\end{array}$$ 
is an equivalence. The left hand-vertical is an equivalence since $K$ is
small. Thus the lower horizontal is an equivalence, 
which is to say that the counit 
$$TEX=\HomR (k,X)\tensorE k \lra X $$ 
is a $K$-equivalence. By proxy-smallness,  this counit map is a $k$-equivalence.
\end{proof}

\begin{examples}
\label{egcell}
(i) If $R$ is a commutative local ring, we saw in \ref{rk:cellular} 
that the $k$-cellular approximation 
of a module $M$ is $\Gamma_IM =K(I) \tensorR M$, so we have
$$TEM \simeq K(I) \tensor_R M.$$

(ii) If $R=C^*(X;k)$ it is not easy to say what the $k$-cellular approximation  is
 in general, but any bounded below module  $M$  is cellular.

(iii) In the equivariant case with $k=F(G_+,R)$ we see that the 
 $k$-cellular approximation  of an $R$-module $M$ is $M \sm EG_+$. Even without
 knowing that $k$ is proxy-small, it is not hard to see that every 
 $G$-free module is effectively constructible. 
\end{examples}

\subsection{Second variant.}

There is a second adjunction  between the derived 
categories of left $\boldR$-modules and  of right $\cE$-modules. 
In the first variant, $k$ played a central role as a left $R$-module
and a left $\cE$-module. In this second variant 
$$\ksharp :=\HombR (k,R)$$
plays a corresponding role: it is a right $R$-module and a right 
$\cE$-module. We have  the adjoint pair 
$$\adjunction{\Ep}{D(\bRmod)}{D(\modcE)}{C}$$
defined by 
$$\Ep (M):=\ksharp \tensorR M \mbox{ and } C(X):=\HomE (\ksharp ,X).$$ 

\begin{remark}
\label{defn:comp}
 If $k$ is small then 
$$\Ep (M) = \HomR (k,R) \tensorR M \simeq \HomR(k, M)=EM, $$
so the two Morita equivalences consider the left and right adjoints of the 
same functor.
\end{remark}

The unit of the adjunction $M \lra C\Ep (M)$ is not
very well behaved, and the functor $C\Ep$ is not even idempotent in general. 
%If $k$ is proxy-small, there is some good behaviour since  
%$M \lra C\Ep (M)$ is completion for 
%$M=R$ and hence for any small $R$-module $M$. 

%\begin{lemma}  
%Provided $\boldk$ is proxy-small, 
% the unit $M \lra C\Ep M$ is $k$-completion 
%and in particular an equivalence if $M$ is cobuilt from $\boldk$. 
%\end{lemma}

\subsection{Complete modules and torsion modules} \label{sec:Rmodfunctors}

Even when we are not interested in the intermediate category 
of $\cE$-modules, several of the composite functors give interesting
endofunctors of the  category 
of $R$-modules.

\begin{lemma}
\label{cellsmashing}
If $k$ is proxy-small then $k$-cellular approximation  is smashing:
$$\cell_kM \simeq (\cell_k R)\tensorR M.$$
\end{lemma}

\begin{proof}
Both modules come with natural maps to $M$. Since $R$ builds $M$, 
$(\cell_kR) \tensorR R$ builds $(\cell_k R) \tensorR M$, so that 
there is a unique map $(\cell_kR) \tensorR M \lra \cell_k M$.

Assuming $k$ is proxy-small, there is a unique map in the reverse 
direction, because $\HomR (\cell_k M, \check{C}(k) \tensorR M)\simeq 0$.
Indeed, this obstruction module cobuilt from $\HomR (K, \check{C}(k) \tensorR M)$, 
and $K^{\#} \tensorR \check{C}(k) \simeq 0$.
\end{proof}

We therefore see by \ref{proxycell} and \ref{cellsmashing} that if $k$ is proxy-small 
$$\cell_k (M)=TE M =T\Ep M.$$
This is the composite of two left adjoints, focusing 
attention on its right adjoint $CE M$, and we note that
$$CE(M)=\HomR (\ksharp, \HomR (k, M))=\HomR (TE R, M).$$
By analogy with Subsection \ref{subsec:loccomp}, we may make the 
following definition.

\begin{defn} 
The {\em completion } of an $R$-module $M$ is the map
$$M \lra \HomR (TER, M)=CEM.$$
We say that $M$ is {\em complete} if the completion 
map is an equivalence. 
\end{defn}

\begin{remark}
By \ref{cellsmashing} we see that completion is idempotent. 
\end{remark}

We  adopt the notation 
$$\Gk M := TE'M$$
and
$$\Lk M := C EM. $$
This is  by analogy with the case of commutative algebra through 
the approach of Part 1 (see Subsection \ref{subsec:loccomp}), where 
$\Gk =\Gamma_I$ is the total 
right derived functor of the $I$-power torsion functor and 
$\Lk =\Lambda^I$ is the total left derived functor of the completion 
functor (see \cite{AJL1,AJL2} for the  context of commutative rings).

It follows from the adjunctions described in Section \ref{sec:Morita} 
that $\Gk$ is left adjoint to $\Lk$ as endofunctors of the category of 
$R$-modules: 
$$\HomR (\Gk M, N)=\HomR (M, \Lk N)$$ 
for $R$-modules $M$ and $N$. Slightly more general is the following observation. 

\begin{lemma} 
\label{Lipman}
For $R$-modules $M$ and $N$ there is a 
natural equivalence
$$\Lk \HomR (M,N) \simeq \HomR (\Gk M, \Gk N).$$
\end{lemma}

\begin{proof} Consider $\HomR (M, \HomR (\Gk R, N))$. \end{proof}

\begin{lemma}
If $k$ is proxy-small, $\Gk$ and $\Lk$ give an adjoint equivalence
$$\cell (R,k) \simeq D(\mathrm{comp}\!-\!\Rmod),$$
where $D(\mathrm{comp}\!-\!\Rmod)$ is the triangulated subcategory 
of $D(R)$ consisting of complete modules.
\end{lemma}

\begin{proof}
We have
$$T\Ep M \simeq T E M \simeq \Gk M \simeq \Gamma_KM, $$
and 
$$CE M \simeq \HomR (\Gk R,M) \simeq \HomR (\Gamma_K R,M),$$
so it suffices to prove the result when $k$ is small. When 
$k$ is small the present adjunction is the composite of two  adjoint pairs of
equivalences. We have seen this for the first variant, and the second
 variant is proved similarly by arguing that the unit and counit
are equivalences. 
\end{proof}

\section{Matlis lifts.}

Matlis lifts are the second novel ingredient in the 
story. They  can be viewed as a language for discussing 
orientations. For a more extensive discussion, see  \cite{DGI}.

\subsection{The definition.}
To motivate the definition, suppose $R$ is a $k$-algebra. There
are then two dualities on $R$-modules we may wish to consider:
the Spanier-Whitehead dual 
$$M^{\#}=\HomR (M,R)$$ 
and the Brown-Comenetz dual 
$$M^{\vee}=\Hom_k(M,k).$$ 
To see how different these usually 
are,  we may conider the case when  $R$ is concentrated in infinitely 
many positive degrees: then $R^{\#}=R$ is still positively graded and 
cyclic, whereas $R^{\vee}$ is negatively graded and not finitely generated.

When $R$ is not a $k$-algebra, the situation is a little more complicated, 
and a guiding example is provided by Pontrjagin duality.
For this we consider the case $R=\Z $ and $k=\Z /p$.
If $T$ is a $\Z /p$-module, we have the vector space dual
$\Hom_{\Z /p}(T,\Z /p)$, but we may wish to extend this
to $\Z$-modules, using
$$\Hom_{\Z /p}(T,\Z /p)=\Hom_{\Z}(T,\Z /p^{\infty});$$
we say $\Z /p^{\infty}$ is a {\em Matlis lift} of
$\Z /p$ from $\Z/p$-modules to $\Z$-modules. 

\begin{defn}
Given a ring map $R \lra k$, and a $k$-module $N$, we say 
that the $R$-module $I=I(N)$ is a Matlis lift of $N$ if
(i) $\HomR (k,I)\simeq N$ as left $k$-modules and (ii)
$I$ is effectively constructible from $k$. 
\end{defn} 

The first condition is the special case $T=k$ of our motivation, 
and the general case follows. Since the first condition only 
sees the $k$-cellularization of $I$, it is natural to require
$I$ is $k$-cellular, and the second condition is a small 
strengthening of this. 

It is worth highlighting the fact
that the existence of a Matlis lift of $N$ provides a right
$\cE$-action on $N$ extending its left $k$-module structure.
Discussion of Matlis lifts can be reformulated in these terms.

We will only need the special case $N=k$ here. If $k$ is proxy-small,
effective constructibility is the same as cellularity so  
\ref{proxycell} shows that Matlis lifts of $k$ correspond to 
right $\cE$-actions on $k$. Given a preferred Matlis lift $I$ of $k$, 
we write  $M^{\vee}=\HomR (M,I)$ since it extends $k$-duality 
of $k$-modules.

\subsection{Examples.} 

 In the particular case of the $k$-module $k$ itself, we may ask how many structures 
it admits as a right $\cE$-module, and for which of these structures a Matlis lift 
exists. In several interesting cases, there is a unique right action of $\cE$ on $k$.

\begin{examples}
\label{egMatlis}
(i) If $R$ is a commutative local ring, the injective hull $I(k)$ is 
a Matlis lift of $k$. Connectivity arguments show there is 
a unique right $\cE$-action on $k$, and hence $I(k)$ is the unique
Matlis lift.

(ii) If $R=C^*(X;k)$, for a connected space $X$, 
 then $\cE \simeq C_*(\Omega X;k)$
provided the Eilenberg-Moore spectral sequence converges.
Right $\cE$-actions on $k$ correspond to actions of 
$\pi_0(C_*(\Omega X;k))=k[\pi_1(X)]$ on $k$, which is to say
group homomorphisms $\pi_1(X) \lra k^{\times}$. The trivial action 
is the one specified by the trivial homomorphism, and the associated
Matlis lift is $C_*(X;k)$.
If $\pi_1(X)$ is a $p$-group and $k=\Fp$, the homomorphism is 
 necessarily trivial, so $C_*(X;k)$ is the unique Matlis lift of $k$. 

(iii) In the equivariant case with $k=F(G_+,R)$, actions of $\cE$
on $k$ correspond to actions of $F(G_+,G_+)$ on $G_+$ extending the
left action of $G_+$, and hence there is only one. 
The Matlis lift of $k$ corresponding to the standard action is $R$ itself. \qqed
\end{examples}

\subsection{Brown-Comenetz duality and Matlis lifts.}
\label{Ip}

The $p$-primary Brown-Comentez dual $I(p)$ of the sphere spectrum 
in $p$-local stable homotopy theory is defined by the condition 
$$[T, I(p)]^n=\Hom_{\Z}(\pi_n(T), \Q /\Z_{(p)}). $$
In fact it is the unique  Matlis lift of $H\Fp$ along  $R=S^0 \lra H\Fp=k$ and
$$I(p)=H\Fp \tensorE H\Fp , $$ 
where  $\cE = \End_{S^0}(H\Fp)$ is the mod $p$ Steenrod algebra ring spectrum. 

%Consequently,  there is a spectral sequence
%$$      \mathrm{Tor}_i^{\cA_p} (\Fp , \Fp) 
%\Rightarrow \pi_* (I(p))$$
%of a form Pontrjagin dual to the ordinary Adams spectral sequence. 

\begin{proof}
To see that $I(p)$ is a Matlis lift, we note that by definition 
$$             \Hom_{S^0} (H\Fp, I(p)) = H\Fp.$$
By connectivity of $\cE$, there is a unique right $\cE$-action
on $H\Fp$. 

Since $H\Fp$ is not proxy-small, we must argue directly 
that $I(p)$ is effectively constructible.
The counit of the adjunction gives a comparison map 
$H\Fp \tensorE H\Fp \lra I(p)$. The homotopy groups of both 
sides are calculated by a spectral sequence with $E^2$-term
$\mathrm{Tor}_i^{\cA_p} (\Fp , \Fp)$, on the left it is a 
bar spectral sequence and on the right a dual Adams spectral 
sequence. In fact we may dualize the usual construction of the 
Adams spectral sequence and observe it also gives a calculation 
of the bar spectral sequence, so the map induces an isomorphism 
of $E^2$-terms.
%%To see that  $ I(p)$ is built from $H\Fp$, note from the 
%% definition that
%% $$\pi_n(I(p))
%% \trichotomy
%% {0 & \mbox{ for } n>0}
%% {\Z/p^{\infty} & \mbox{ for } n=0}
%% {\mbox{a finite $p$-group} & \mbox{ for } n <0.}$$
%% By killing homotopy groups we obtain cofibre sequences $$I(p)[-n,0] \lra I(p) \lra I(p)(-\infty, -n-1]$$ and hence we obtain an equivalence $$I(p)\simeq\colim_n I(p)[-n,0]$$ showing that $I(p)$ is built from $H \Fp$. 

%% It is not hard to see that there is a unique action on $H\Fp$ and hence
%% at most one Matlis lift.
%% $H\Fp \smE H \Fp$: the conditions on $R=S^0$ and $k=H\Fp$ are obvious, 
%% and the condition on $N\tensorE k$ follows from the spectral sequence 
%%$$\Tor^{\cA_p}_{*,*}(\Fp , \Fp ) \Rightarrow \pi_*(H\Fp \smE H\Fp ).$$
\end{proof}

\subsection{Local duality and the dualizing complex.} \label{sec:localduality}

A dualizing complex $D$ for  a  commutative Noetherian local ring $R$ 
is characterized by the three requirements
(i) $\HomR (D,D)=R$, (ii) $D$ has finite injective dimension 
and (iii) $H_*(D)$ is a finitely generated $R$-module. 
There is a useful summary of properties in \cite{af2}.

Dualizing complexes exist under very 
weak hypotheses, certainly if the ring is a quotient of 
a complete local ring. 
We only use the fact that for a local ring $R$ the dualizing complex $D$ 
has the property that $\Gk D =\I (k)$ (up to suspension). From 
the fact that $\Gk$ and $\Lk$ give inverse equivalences, 
it is therefore natural to make the following definition.

\begin{defn}
If there is a Matlis lift $\I (k)$ of $k$, 
we define the {\em dualizing complex} by 
$$D=\Lk \I (k).$$
\end{defn}

\begin{remark} To justify calling $D$ the dualizing complex, 
we should establish a natural equivalence
$$R \lra \HomR (D,D)=\HomR (\I (k), \I (k)),$$
in some generality, at least when $R$ is complete. Here
are two examples from local algebra.

(a) Suppose $R$ is a $k$-algebra. In this case 
$\HomR (\I (k), \I (k))=\Homk (\I (k),k)$. 
If $R$ is $k$-reflexive (for example if $k$ is a field and 
$H_*(R)$ is bounded below and a finite dimensional vector space 
in each degree) we also have  $\I (k)=\Homk (R,k)$, and the
required equivalence follows.

(b) Suppose $R$ is a complete, commutative local ring.
Suppose also that there is a map $Q \lra R$ making $R$ into a 
finite $Q$-module and with $Q$ complete and  Gorenstein. 

In that case $\Gamma_{IQ}=\Gamma_I$ on $R$-modules, and $\Gamma_IQ=\I_k^Q(k)$ is the Matlis lift of $k$ to a 
$Q$-module. Hence
$$Q=\Hom_Q(\Gamma_I Q, \Gamma_I Q)=\Hom_Q(\I_k^Q(k),\I_k^Q(k)).$$ We now claim 
$$\I_k^R(k):=\Hom_Q(R,\I_k^Q(k))$$
is a Matlis lift of $k$ to an $R$-module. We may then calculate $$ R=Q \tensor_Q R=\Hom_Q(\I_k^Q(k),\I_k^Q(k))\tensor_Q R= \Hom_R(\I_k^R(k),\I_k^R(k)).\qqed$$
\end{remark}

We may connect the dualizing complex with local duality by a 
tautology.

\begin{cor}
\label{locduality}
(Local duality)
If  $\I (k)$ is a Matlis lift of $k$ there
is an equivalence
$$\HomR (\Gk R, \I (k))=\Lk D.\qqed$$
\end{cor}

%\begin{proof}
%Taking $M=R$ and $N=D$ in \ref{Lipman} we 
%obtain $\Lambda^k\HomR (R,D) \simeq \HomR (\Gamma_k R, \Gamma_k D)$. 
%Since $\I (k)$ is cellular, $\Gamma_k D =\Gamma_k\Lambda^k \I (k) \simeq \I (k)$.
%\end{proof}

\section{The Gorenstein condition.}
A commutative Noetherian local ring $R$ is Gorenstein if and only if $\Ext_{R}^*(k,R)$ 
is  one dimensional as a $k$-vector space \cite[18.1]{Mat}. In the derived category we 
can restate this as saying that the homology of the (right derived) Hom complex 
$\Hom_{R}(k,R)$ is equivalent to a suspension of $k$ (cf \cite{fht}). 
This suggests the definition for ring spectra.

\begin{defn} \cite{DGI}
We say that $\boldR\lra \boldk$ is {\em Gorenstein} if there is an equivalence of 
$\boldR$-modules 
$\Hom_\boldR(\boldk,\boldR) \simeq \Sigma^a\boldk$ for some integer $a$. 
\end{defn}

As in commutative algebra, Gorenstein rings are ubiquitous.

\begin{examples}
(i) The ring spectrum associated to a commutative local ring $R$
is Gorenstein if and only if the local ring is  Gorenstein in the 
conventional sense, as we discussed in our motivation.

(ii) The ring spectrum $R=C^*(X;k)$ is Gorenstein in two cases. 
 First there is the more familiar case when $X$ is a compact connected
 manifold orientable over $k$.  In fact $R$ is also Gorenstein 
if $X$ is not orientable, provided $2^n$ acts as zero on $k$ for some $n$.
To see this, let $\tR =C^*(\tX;k)$ denote the cochain complex of the
orientable double cover $\tX$ of $X$, and let $C$ be the group of order 2, acting 
via covering transformations. Notice that the convergence of the 
Rothenberg-Steenrod and Eilenberg-Moore spectral sequences give
 $\Hom_{kC}(k,\tR) \simeq R$ and 
$\tR \tensorR k \simeq kC$. Now we may calculate
$$\begin{array}{rcl}
\HomR (k,R) & \simeq & \HomR (k, \Hom_{kC}(k,\tR))\\
 & \simeq & \HomR (k \tensor_{kC} k, \tR)\\
 & \simeq & \HomtR (\tR \tensorR k \tensor_{kC} k, \tR)\\
& \simeq & \HomtR (kC \tensor_{kC} k, \tR)\\
& \simeq & \HomtR ( k, \tR).
\end{array}$$

 Secondly,  $R$ is Gorenstein when $X=BG$ for any finite group $G$.
 We will give the proof when $G$ is a finite $p$-group in 
 Subsection \ref{sec:MoritaGorenstein} below. In fact this example
 extends to many cases when $G$ is a compact Lie group. To explain, 
 we need to consider  the adjoint action of 
$\pi_1(BG)$ on $H^d_c(T_eG;k)\cong k$, where $G$ is of dimension $d$.  The ring 
 $R$ is Gorenstein if  the action is trivial or if $2^n$ acts as 0 on $k$.

(iii) The group ring $kG$ of a finite group $G$ is Gorenstein as a ring spectrum.
More generally, if $R$ is an equivariant commutative ring spectrum and
$k=F(G_+,R)$ represents the nonequivariant theory, then 
$$\HomR (k,R)=\HomR(F(G_+,R),R)\simeq G_+\sm \HomR (R,R) \simeq k,$$
so that $R$ is Gorenstein.
\end{examples}

\subsection{Orientability.}
For simplicity we suppose for the rest of the section 
that all cellular objects are effectively constructible, for 
example if $k$ is proxy-small. 

If $\boldR$ is Gorenstein, $\boldk$ acquires new structure: 
that of a right $\cE$-module (not to be confused with its natural structure 
as a {\em left} $\cE$-module), and $\Gamma_k R$ is a Matlis lift of $k$.
We want to say that $\boldR$ is {\em orientable} if the new right action is `trivial'.
In many examples there is a natural candidate for an action to be called 
trivial, and a corresponding Matlis lift $I$ of $k$.

\begin{defn}
A Gorenstein ring spectrum $\boldR$ is {\em orientable} if 
$$\HombR(\boldk,\boldR) \simeq \Sigma^a \HombR (\boldk,I)$$
 as right $\cE$-modules. 
\end{defn}

\begin{examples}
(i) If $R$ is a commutative local ring, we remarked in 
Example \ref{egMatlis} that there is a unique right $\cE$-module 
structure on $\boldk$, and hence a unique 
Matlis lift of $k$. Thus every Gorenstein commutative ring is orientable 
as a ring spectrum. The notion of orientability is redundant in 
classical commutative algebra. 

(ii) If $R=C^*(X;k)$, we have already remarked that 
the Matlis lift corresponding to the trivial action is $I=C_*(X;k)$.

If $X$ is a manifold for which $R=C^*(X;k)$ is Gorenstein, the associated
 action is the one arising from the monodromy action. Thus 
$R$ is orientable if and only if  the manifold $X$ is orientable over $k$.
For example if $X$ is not orientable and $k=\Z /4$, the ring spectrum $R$ is
 Gorenstein but not orientable.

Similarly, if $X=BG$, the Gorenstein action is the adjoint action of 
$\pi_1(BG)$ on $H^d_c(T_eG;k)$, where $G$ is of dimension $d$. 
The case $G=O(2)$ gives an example where this is non-trivial, since
the reflections  act as  $-1$. 
The ring  $k=\Z /4$ is still proxy-small so $R$ is Gorenstein 
but not orientable.

(iii) In the equivariant example with $k=F(G_+,R)$, 
there is an obvious action of $F(G_+,G_+)$
on $G_+$, giving an action of $\cE$ on $k$, which we use to give the Matlis
 lift. The Gorenstein action is its opposite; since the two are isomorphic,
this example is orientably Gorenstein. 
\end{examples}

\subsection{The local cohomology theorem.}
We saw in \ref{rk:cellular} that the  stable Koszul complex provides a construction of
$\boldk$-cellular approximation in many contexts. In this case,  its homotopy is calculated 
using local cohomology: there is a spectral sequence 
$$
H^*_{\fm}(\pi_*(\boldR))\Rightarrow \pi_*(\Gamma_k \boldR).
$$

 We can then deduce a valuable duality property from the Gorenstein condition. 
Indeed, if $\boldR \lra \boldk$ is Gorenstein and orientable, we have the equivalences 
$$ E\Gamma_k \boldR =\HombR (\boldk,\Gamma_k \boldR) \simeq 
\HombR(\boldk,\boldR) \simeq \Sigma^a \HombR (\boldk, I) =E\Sigma^aI $$ 
of right $\cE$-modules. By \ref{proxycell},   we conclude
$$\Gamma_k \boldR \simeq TE\Gamma_k \boldR \simeq TE \Sigma^aI \simeq \Sigma^a I.$$ 
For example if $\boldR$ is a $\boldk$-algebra for a field $k$, 
we can take $I =\Homk (\boldR,\boldk)$ provided it is $k$-cellular, 
and conclude that there is a spectral sequence 
$$H^*_{\fm}(\pi_*(\boldR))\Rightarrow H^*(\Homk (\boldR,\boldk))=\Homk (\pi_*\boldR,k).$$ 
In particular {\em if $\pi_*(\boldR)$ is Cohen-Macaulay}, this spectral sequence 
collapses to show {\em it is also Gorenstein}. In fact one may apply Grothendieck's dual 
localization process to this spectral sequence and hence conclude that 
{\em whatever its depth,  $\pi_*(\boldR)$ is generically Gorenstein} \cite{ringlct}.
This dual localization process is lifted to the level of module spectra in 
\cite{kappaI}.

\begin{examples}
(i) When $R$ is a commutative local ring, we just obtain Gorenstein duality, 
stating that a Gorenstein ring is Cohen-Macaulay and $H^r_{\fm}(R)=I(k)$.

(ii) If $R=C^*(X;k)$ for an orientable manifold $X$ 
we deduce the Poincar\'e duality statement
$$H^*(X;k)=H^0_{\fm}(H^*(X;k))\cong H_*(X;k).$$

A much more interesting example is that of $C^*(BG)$ for a compact Lie group $G$. 
Since $\pi_*(C^*(BG))=H^*(BG)$,  when $G$ is finite there is a spectral sequence 
$$H^*_{\fm}(H^*(BG))\Rightarrow H_*(BG), $$ 
showing that the group cohomology ring $H^*(BG)$ has very special properties, 
such as being generically Gorenstein. 

(iii) For the group ring, we conclude $kG =H^0_{\fm}(kG)=kG^{\vee}$, which is
to say that $kG$ is a Frobenius  algebra.

 For the general equivariant example, with $k \simeq F(G_+,R)$, 
 we have $\Gamma_kR\simeq R \sm EG_+$.
However, in this case $\Homk (R,k) \simeq R$, which is not $k$-cellular:
the statement $\Gamma_kR \simeq \Gamma_k (R^{\vee})$ is a tautology. The only 
content in this case comes from the work of Section \ref{sec:compthms1} of Part 1, 
where we say $R \sm EG_+ \simeq \Gamma_k R \simeq \Gamma_I R$, converting the 
geometric cellular approximation  to an algebraic one.
\end{examples}

\subsection{The Gorenstein condition and the dualizing complex.}

Another characterization of the Gorenstein condition for 
complete local rings is condition on $R$ that  $R \simeq \Sigma^a D$
for some $a$, where $D$ is the dualizing complex. 
Combining this with the tautologous local duality statement \ref{locduality}, 
we see that this version of the Gorenstein condition implies
$$\HomR (\Gk R, \I (k))= \Sigma^aR.$$

\begin{example}
 The example of Mahowald-Rezk \cite{MR} fits well here. 
They define the $n$th dualizing complex by 
$$W_nS^0=IC_n^fS^0,  $$
where $I$ denotes Brown-Comenetz duality, and $C_n^f$ is the 
finite counterpart of the $n$th acyclization discussed in Section 
\ref{sec:bordism}, namely the $F(n+1)$-cellular
approximation functor, for a finite type $n$-spectrum $F(n+1)$.
Thus it is natural to write $C_n^fS^0 =\Gamma_{F(n+1)}S^0$, and hence 
its Brown-Comenetz dual is the dualizing complex. 

Mahowald and Rezk consider the class of spectra for which their cohomology 
is finitely presented over the Steenrod algebra. They prove that all such 
type $n$ fp spectra  are reflexive. In 
fact they go on to form $WX$ as a direct limit of $W_nX$ and show that for many 
interesting fp ring spectra (such as $X=ku,ko, eo_2, BP \langle n \rangle$)
we have $WX \simeq \Sigma^dX$ for a suitable $d$. This says that the dualizing 
complex is a suspension of the ring, and we have just observed this is a form
of the Gorenstein condition.
\end{example}

\subsection{Morita invariance of the Gorenstein condition.} 
\label{sec:MoritaGorenstein}

In this final section we show that the Gorenstein condition is Morita invariant 
in many useful cases, provided $R$ is a $k$-algebra. This allows us to 
deduce striking consequences from well-known examples of Gorenstein rings. 
For instance we can deduce the local cohomology theorem for finite $p$-groups 
from the fact that $kG$ is a Frobenius algebra.

\begin{thm}
\label{REGorenstein}
Suppose $R$ is a $k$-algebra, and $R^{\vee}$ is $k$-cellular and
finally that  $\cE$ and $R$ are Matlis reflexive.
Then
$$\HomE (k , \cE )\cong \HomR (k,R),$$
and hence
$$\cE \mbox{ is Gorenstein } \Longleftrightarrow  R \mbox{ is Gorenstein. }$$ 
\end{thm}

\begin{proof} We use the fact that
$E(R^{\vee})=\kvee $, so that we have 
$$R^{\vee}=TE(R^{\vee})=T\kvee =\kvee \tensorE k, $$
and the fact that $\cE =\HomR (k, k)$ dualizes to give 
$$\cE^{\vee}=k  \tensorR \kvee . $$ 
%This is an isomorphism of $\cE$-bimodules. 
%However the left $\cE$-module structure comes from $k $ and the right 
%module structure from $\kvee$.

Next, note that   the expression $k \tensorR \kvee \tensorE k$
makes sense, where the right $\cE$-module structure on the first two factors comes from $\kvee$.  The key equality in the proof is simply the 
associativity isomorphism
$$\cE^{\vee} \tensorE k=
   k \tensorR \kvee \tensorE k=k \tensorR R^{\vee} .$$

Now we make the following calculation,

$$\begin{array}{rcl}
\HomE(k,\cE )   &\simeq&\HomE (k , (\cE^{\vee})^{\vee})\\
                &\simeq&\Homk ( \cE^{\vee}\tensorE k ,k)\\
                &\simeq&\Homk (k \tensorR \kvee \tensorE k,k)\\
                &\simeq&\Homk (k \tensorR R^{\vee}  ,k)\\
                &\simeq&\HomR (k ,(R^{\vee})^{\vee})\\
                &\simeq&\HomR (k ,R)
\end{array}$$
\end{proof}

\begin{cor}
If $G$ is a finite $p$-group then $kG$ is Gorenstein and hence
$C^*(BG)$ is Gorenstein.\qqed
\end{cor}

\end{document}